\documentclass[preprint,3p,12pt]{elsarticle}
\usepackage{mathrsfs}
\usepackage{amsmath}
\usepackage{stmaryrd}
\usepackage{bbding}
\usepackage{dcolumn}
\usepackage{graphicx}
\usepackage{amsfonts}
\usepackage{amssymb}
\usepackage{psfrag}
\usepackage{wrapfig}
\usepackage{subfigure}
\usepackage{makeidx}
\usepackage{bm}
\usepackage{epsf}
\usepackage{epsfig}
\usepackage{setspace}
\usepackage{graphicx}
\usepackage{epstopdf}
\usepackage{psfrag}
\usepackage{subfigure}
\usepackage{color}

\begin{document}

\title{A sixth-order weighted essentially non-oscillatory scheme for hyperbolic conservation laws}

\author[caep]{Fengxiang Zhao}
\ead{kobezhao@126.com}

\author[iapcm]{Liang Pan}
\ead{panliangjlu@sina.com}

\author[iapcm]{Zheng Li}
\ead{zheng\_li@iapcm.ac.cn}

\author[iapcm]{Shuanghu Wang\corref{cor}}
\ead{wang\_shuanghu@iapcm.ac.cn}

\address[caep]{The Graduate School of China Academy of Engineering Physics, Beijing, China}
\address[iapcm]{Institute of Applied Physics and Computational Mathematics, Beijing, China}
\cortext[cor]{Corresponding author}

\begin{abstract}
In this paper, A new sixth-order weighted essentially non-oscillatory (WENO) scheme, refered as the WENO-6, is proposed in the finite volume framework for the hyperbolic conservation laws. Instead of selecting one stencil for each cell in the classical WENO scheme \cite{WENO-JS}, two independent stencils are used for two ends of the considering cell in the current approach. Meanwhile, the stencils, which are used for the reconstruction of variables at both sides of interface, are symmetrical. Compared with the classical WENO scheme \cite{WENO-JS}, the current WENO scheme achieves one order of improvement in accuracy with the same stencil. The reconstruction procedure is defined by a convex combination of reconstructed values at cell interface, which are constructed from two quadratic and two cubic polynomials. The essentially non-oscillatory property is achieved by the similar weighting methodology as the classical WENO scheme. A variety of numerical examples are presented to validate the accuracy and robustness of the current scheme.
\end{abstract}


\maketitle

\section{Introduction}
In past decades, there have been tremendous efforts on designing high-order accurate numerical schemes for compressible fluid flows and great success has been achieved. High-order accurate numerical schemes were pioneered by Lax and Wendroff \cite{Lax-Wendroff}, and extended into the version of high resolution methods by van Leer \cite{Van-Leer}, Harten \cite{Harten} et al. and other higher order versions, such as essentially non-oscillatory scheme (ENO) \cite{ENO1, ENO2}, weighted essentially non-oscillatory scheme (WENO) \cite{WENO-JS, WENO}, Hermite weighted essentially non-oscillatory scheme (HWENO) \cite{HWENO1, HWENO2}, and discontinuous Galerkin scheme (DG) \cite{DG2, DG3, DG1}, etc.

The ENO and WENO schemes have been successfully applied for the compressible flows with strong shocks, contact discontinuities and complicated smooth structures. The ENO schemes was first introduced by Harten et al. \cite{ENO1} in the form of cell averaged variables. The key idea of ENO schemes is to use the "smoothest" stencil among several candidates to approximate the fluxes at cell interfaces to achieve high order accuracy and avoid spurious oscillations near discontinuities. Later, the flux version of ENO schemes \cite{ENO2} was introduced with TVD Runge-Kutta temporal discretization. However, the ENO scheme is not effective in terms of selecting only one stencil to approximate the fluxes at cell interface, and such a adaption of stencils is not necessary in smooth regions. To overcome these drawbacks while keeping the robustness and high order accuracy of ENO scheme, the WENO scheme was first introduced in \cite{WENO}. Instead of approximating the pointwise value of the solution using only one of the candidate stencils, a convex combination of all the candidate stencils was used. Each candidate stencil is assigned a weight which determines the contribution of this stencil to the final approximation of pointwise value. The weights can be defined in such a way that it approaches certain optimal weights to achieve a higher order of accuracy in smooth regions, and the stencils which contain the discontinuities are assigned a nearly zero weight in the regions near discontinuities. A higher order of accuracy is obtained by emulating upstream scheme with the optimal weights away from the discontinuities, and the essentially non-oscillatory property is achieved near discontinuities. The flux of WENO scheme is smoother than that of the ENO scheme, and the smoothness enables us to prove convergence of WENO scheme for smooth solutions using Strang's technique \cite{WENO-JS}.

However, the optimal order of WENO scheme was not attained in \cite{WENO}, i.e. $(2r-1)$-th order with $r$-th order ENO scheme. A more detailed error analysis for WENO scheme was carried out, and a new WENO scheme (WENO-JS) including smoothness indicators and nonlinear weights was constructed in \cite{WENO-JS}. It was made possible to generalize the scheme to the fifth-order of accuracy. Later, very high order WENO schemes were developed as well \cite{very-high-WENO}. However, the WENO-JS scheme may lose accuracy if the solution contains local smooth extrema. A new WENO scheme (WENO-M) was developed to overcome this problem by modifying the nonlinear weights by a mapping procedure. However, the proposed mapping procedure is revealed to be computationally expensive. With a different weighting formulation, another version of the fifth-order WENO scheme (WENO-Z) was introduced in \cite{WENO-Z}, in which a global higher order reference value was used for the smoothness indicator. As the improvement of WENO-Z scheme, WENO-Z+ scheme was also developed. The WENO-M and WENO-Z and WENO-Z+ schemes turned out to be less dissipative than the classical WENO-JS scheme near smooth extrema. In the solution with high-frequency waves, they achieve noticeably higher amplitudes than WENO-JS scheme in a coarse grid.

In this paper, a new six-order WENO scheme in finite volume framework is developed by a nonlinear convex combination approach with all corresponding polynomials to obtain high-order point-wise values at cell interface.
As a complementary version of the classical WENO schemes based on cell averages, a new reconstruction procedure is proposed. At two ends of each cell, two independent stencils are selected to construct the interface values separately. The two independent stencils, which are used for the reconstruction of variables at left and right sides of cell interface, are symmetrical.  The reconstruction procedure is defined by a convex combination of reconstructed values at cell interface, which are constructed from two quadratic and two cubic polynomials. The essentially non-oscillatory property is preserved by the similar weighting methodology as the classical WENO schemes proposed in \cite{WENO-JS}. Compared with the classical WENO schemes, the current WENO scheme can achieve one order of improvement in accuracy and better resolution power with the same stencil, while preserving a good robustness.

This paper is organized as follows. In section 2, the classical WENO scheme in finite volume framework are briefly reviewed. The new WENO scheme is introduced in section 3, where a detailed discussion is also given. Section 4 includes numerical examples to validate the current algorithm. The last section is drawing the conclusion.

\section{Finite volume type WENO scheme}

\subsection{Finite volume methods}
We consider the following hyperbolic conservation law
\begin{equation}\label{hyperbolic}
\frac{\partial W}{\partial t}+\frac{\partial F(W)}{\partial x}=0,
\end{equation}
with the initial condition
\begin{equation*}
W(x,0)=W_{0}(x).
\end{equation*}
Integrating Eq.\eqref{hyperbolic} over cell $I_i=[x_{i-1/2}, x_{i+1/2}]$, the semi-discretized form of finite volume scheme can be written as
\begin{align}\label{ode}
\frac{\text{d}W_{i}}{\text{d}t}=\mathcal{L}_i(W(x))=-\frac{1}{\Delta x}[F(W(x_{i+1/2},t))-F(W(x_{i-1/2},t))],
\end{align}
where $W_{i}$ is the cell averaged value, $\Delta x$ is the cell size, and $F(W(x_{j+1/2},t))$ is the flux at cell interface $x=x_{i+1/2}$, which can be approximated by numerical flux $F_{i+1/2}$ as follows
\begin{align*}
F(W(x_{i+1/2},t))\approx F_{i+1/2}=F(W_{i+1/2}^{l},W_{i+1/2}^{r}),
\end{align*}
where $W_{i+1/2}^{l}$ and $W_{i+1/2}^{r}$ are the reconstructed values at both sides of the cell interface. To fully discretize Eq.\eqref{hyperbolic}, the approximate Riemann solvers can be used \cite{Riemann-appro} for the numerical flux, the classical third-order TVD Runge-kutta scheme\cite{TVD-rk} and two-stage fourth-order scheme \cite{fourth-GRP} can be used for temporal discretization. The spatial discretization is the main theme of this paper, which will be given in the following sections.

\subsection{the classical WENO scheme}
Before the new WENO scheme is introduced, we will briefly review the classical WENO reconstruction \cite{WENO} in this section. Assume that $W(x)$ are the variables which need to be reconstructed, $W_i$ are the cell averaged values, and $W_i^r, W_i^l$ are the two values obtained by the reconstruction at two ends of the $i$-th cell. The fifth-order WENO reconstruction is given as follows
\begin{align*}
W_i^r=\sum_{k=0}^2\delta_k^r w_{k}^r,~~
W_i^l=\sum_{k=0}^2\delta_k^l w_{k}^l,
\end{align*}
where all quantities involved are taken as
\begin{align*}
w_{0}^r= \frac{1}{3}W_{i}+\frac{5}{6}W_{i+1}-\frac{1}{6}W_{i+2},&~~
w_{0}^l=\frac{11}{6}W_{i}-\frac{7}{6}W_{i+1}+\frac{1}{3}W_{i+2},\\
w_{1}^r=-\frac{1}{6}W_{i-1}+\frac{5}{6}W_{i}+\frac{1}{3}W_{i+1},&~~
w_{1}^l=\frac{1}{3}W_{i-1}+\frac{5}{6}W_{i}-\frac{1}{6}W_{i+1},\\
w_{2}^r=\frac{1}{3}W_{i-2}-\frac{7}{6}W_{i-1}+\frac{11}{6}W_{i},&~~
w_{2}^l=-\frac{1}{6}W_{i-2}+\frac{5}{6}W_{i-1}+\frac{1}{3}W_{i},
\end{align*}
and $\delta_k^r$ and $\delta_k^l$ are the nonlinear weights. In order to deal with the discontinuity, the local smoothness indicator $\beta_k$ is introduced. For the fifth-order reconstruction, this definition yields
\begin{align*}
\beta_0&=\frac{13}{12}(W_{i}+2W_{i+1}+W_{i+2})^2 +\frac{1}{4}(3W_{i}-4W_{i+1}+W_{i+2})^2,\\
\beta_1&=\frac{13}{12}(W_{i-1}+2W_{i}+W_{i+1})^2 +\frac{1}{4}(W_{i-1}-W_{i+1})^2,\\
\beta_2&=\frac{13}{12}(W_{i-2}+2W_{i-1}+W_{i})^2 +\frac{1}{4}(W_{i-2}-4W_{i-1}+3W_{i})^2,
\end{align*}
The most widely used is the WENO-JS non-linear weights \cite{WENO-JS}, which can be
written as follows
\begin{align*}
\displaystyle \delta_k^{p,JS}=\frac{\alpha_k^{p,JS}}{\sum_{p=0}^2 \alpha_p^{p,JS}},~ \alpha_k^{p,JS}=\frac{d_k^{p}}{(\beta_k+\epsilon)^2},
\end{align*}
where $p=l, r$, and
\begin{align*}
\displaystyle \delta_0^l=\delta_2^r=\frac{3}{10},
\delta_1^l=\delta_1^r=\frac{3}{5}, \delta_2^l=\delta_0^r=\frac{1}{10},
\epsilon=10^{-6},
\end{align*}
and $\beta_k$ is the smooth indicator. In order to achieve the better performance of WENO scheme near smooth extrema, WENO-Z \cite{WENO-Z} reconstruction was developed. The only difference is the nonlinear weights, and the nonlinear weights for WENO-Z scheme are written as
\begin{align*}
\delta_k^{p,Z}=\frac{\alpha_k^{p,Z}}{\sum_{k=0}^2 \alpha_k^{p,Z}},~
\alpha_k^{p,Z}=d_k^{p}\Big[1+\Big(\frac{\tau}{\epsilon+\beta_k}\Big)\Big],
\end{align*}
where $\tau=|\beta_0-\beta_2|$ is used for the fifth-order reconstruction.

\section{The recursive WENO scheme}
\subsection{Linear reconstruction}
In this section, a six-order recursive WENO scheme will be presented, which shares the identical stencil to the classical fifth-order WENO scheme. For the cell interface $x=x_{i+1/2}$, a unique symmetric stencil $S_{i+1/2}=\{I_{i-2},I_{i-1},I_{i},I_{i+1},I_{i+2},I_{i+3}\}$ can be selected. On the stencil, an optimal fifth degree polynomial approximation
\begin{align}\label{polynomial}
W^{opt}(x)\equiv\sum_{k=0}^5 a_{k}x^{k},
\end{align}
can be uniquely determined by the following conditions
\begin{align*}
\frac{1}{\Delta x}\int_{I_{k}}W(x)\text{d}x=W_k, k=i-2,...,i+3,
\end{align*}
The coefficients $a_{k}, k=0,...,5$ can be found by solving the linear system which is derived by the equation above. Substituting $a_{k}$ into Eq.\eqref{polynomial}, $W^{opt}_{i+1/2}$ at the cell interface $x=x_{i+1/2}$ can be written as
\begin{align}\label{central}
W^{opt}_{i+1/2}=\frac{1}{60}(W_{i-2}-8W_{i-1}+37W_{i}+37W_{i+1}-8W_{i+2}+W_{i+3}).
\end{align}

\begin{figure}[!htb]
\centering
\includegraphics[width=0.8\textwidth]{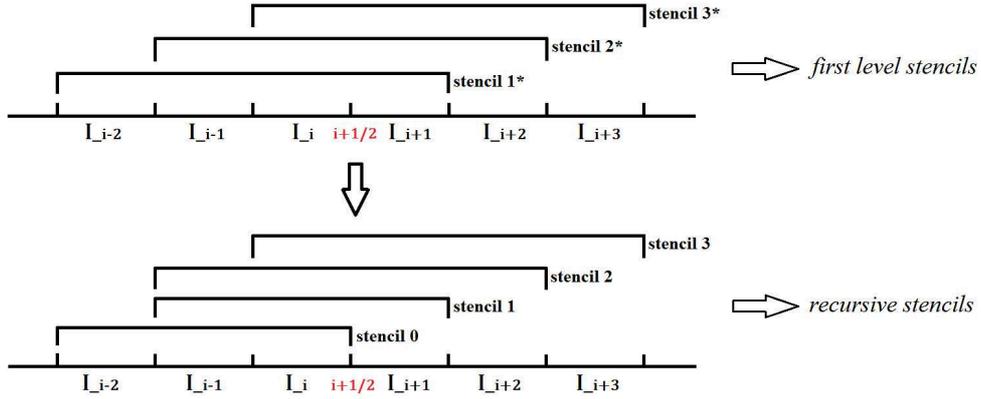}
\caption{\label{stencils} Stencils for the left side of interface $x_{i+1/2}$ in cell $I_{i}$ in sixth-order recursive WENO.}
\end{figure}

In order to achieve high order and non-oscillatory reconstruction, a convex combination of $W_{i+1/2}^{l,r}$ needs to be constructed by the candidate polynomials at both sides of interface. For the reconstruction at the left side of interface $x=x_{i+1/2}$ corresponding to the cell $I_{i}$, the first level stencil and so-called recursive stencil are shown in Fig.\ref{stencils}. In the first level stencil, three sub-stencils containing four cell averaged values are considered, which are denoted by $S^{\star}_{1}, S^{\star}_{2}$ and $S^{\star}_{3}$. In order to deal with discontinuities without oscillation, all possible situations of discontinuities need to be considered in the division of sub-stencils. However, when the discontinuity is at the interface $x=x_{i+1/2}$, all the first level sub-stencils $S^{\star}_{1}, S^{\star}_{2}$ and $S^{\star}_{3}$ are across the discontinuity. $S^{\star}_{1}$ needs to be devided into two "smaller" ones denoted by $S_{0}, S_{1}$, and the $S_{0}$ is the dominant one in dealing with such a situation. Thus, the recursive sub-stencils for sixth-order recursive WENO scheme are given, which are denoted by $S_{0}, S_{1},S_{2}$ and $S_{3}$ as shown in Fig.\ref{stencils}. The corresponding candidate polynomials are given as
\begin{align*}
S_{0,i+1/2}^{l}&=\{I_{i-2},I_{i-1},I_{i}\} \leftrightarrow w^{l}_{0}(x),\\
S_{1,i+1/2}^{l}&=\{I_{i-1},I_{i},I_{i+1}\} \leftrightarrow w^{l}_{1}(x),\\
S_{2,i+1/2}^{l}&=\{I_{i-1},I_{i},I_{i+1},I_{i+2}\} \leftrightarrow w^{l}_{2}(x),\\
S_{3,i+1/2}^{l}&=\{I_{i},I_{i+1},I_{i+2},I_{i+3}\} \leftrightarrow w^{l}_{3}(x),
\end{align*}
where $w^{l}_{0}(x)$ and $w^{l}_{1}(x)$ are quadratic polynomials, and $w^{l}_{2}(x)$ and $w^{l}_{3}(x)$ are cubic polynomials, which can be determined uniquely by
\begin{align*}
\frac{1}{\Delta x}\int_{I_{k}}w^{l}_{k}(x)\text{d}x=W_k, I_{k}\in S_{k,i+1/2}^{l}, k=0,...,3.
\end{align*}
According to the equation above, the polynomials $w^{l}_{k}(x), k=0,...,3$ can be fully determined, and the point values at the left side of interface $x=x_{i+1/2}$ can be written as
\begin{equation*}
\begin{aligned}
w_{0,i+1/2}^{l}&=\frac{1}{6}(2W_{i-2}-7W_{i-1}+11W_{i}), \\
w_{1,i+1/2}^{l}&=\frac{1}{6}(-W_{i-1}+5W_{i}+2W_{i+1}), \\
w_{2,i+1/2}^{l}&=\frac{1}{12}(-W_{i-1}+7W_{i}+7W_{i+1}-W_{i+2}), \\
w_{3,i+1/2}^{l}&=\frac{1}{12}(3W_{i}+13W_{i+1}-5W_{i+2}+W_{i+3}).
\end{aligned}
\end{equation*}
With the these point values, a convex combination for $W_{i+1/2}$ can be derived as follows
\begin{align}\label{left}
W_{i+1/2}^l=\sum_{k=0}^{3}d_{k}^{l}w_{k,i+1/2}^{l},
\end{align}
Comparing the coefficients of Eq.\eqref{central} with that of Eq.\eqref{left}, the linear weights for the left side can be written as
\begin{align*}
d_{0}^{l}=\frac{1}{20},d_{1}^{l}=\frac{3}{20},d_{2}^{l}=\frac{3}{5},d_{3}^{l}=\frac{1}{5}.
\end{align*}
Similarly, the sub-stencils and candidate polynomials for the reconstruction of right side of interface $x=x_{i+1/2}$ corresponding to the cell $I_{i+1}$ are given as
\begin{align*}
S_{0,i+1/2}^{r}&=\{I_{i+1},I_{i+2},I_{i+3}\} \leftrightarrow w^{r}_{0}(x),\\
S_{1,i+1/2}^{r}&=\{I_{i},I_{i+1},I_{i+2}\} \leftrightarrow w^{r}_{1}(x),\\
S_{2,i+1/2}^{r}&=\{I_{i-1},I_{i},I_{i+1},I_{i+2}\} \leftrightarrow w^{r}_{2}(x),\\
S_{3,i+1/2}^{r}&=\{I_{i-2},I_{i-1},I_{i},I_{i+1}\} \leftrightarrow w^{r}_{3}(x),
\end{align*}
and the linear convex combination can be derived as follows
\begin{align}\label{right}
W_{i+1/2}^r=\sum_{k=0}^{3}d_{k}^{r}w_{k,i+1/2}^{r},
\end{align}
where these point values in the right cell can be written as
\begin{align*}
w_{0,i+1/2}^{r}&=\frac{1}{6}(11W_{i+1}-7W_{i+2}+2W_{i+3}), \\
w_{1,i+1/2}^{r}&=\frac{1}{6}(2W_{i}+5W_{i+1}-W_{i+2}), \\
w_{2,i+1/2}^{r}&=\frac{1}{12}(-W_{i-1}+7W_{i}+7W_{i+1}-W_{i+2}), \\
w_{3,i+1/2}^{r}&=\frac{1}{12}(W_{i-2}-5W_{i-1}+13W_{i}+3W_{i+1}),
\end{align*}
and the linear weights for the right side are given as
\begin{align*}
d_{0}^{r}=\frac{1}{20},d_{1}^{r}=\frac{3}{20},d_{2}^{r}=\frac{3}{5},d_{3}^{r}=\frac{1}{5}.
\end{align*}

\subsection{Nonlinear weights}
With the linear weights, Eq.\eqref{left} and Eq.\eqref{right} are unable to deal with discontinuity without spurious oscillations. In order to overcome this problem, the nonlinear weights are introduced and Eq.\eqref{left} and Eq.\eqref{right} are modified as
\begin{align}\label{right-left}
W_{i+1/2}^l=\sum_{k=0}^{3}\delta_{k}^{l}w_{k,i+1/2}^{l},~ W_{i+1/2}^r=\sum_{k=0}^{3}\delta_{k}^{r}w_{k,i+1/2}^{r},
\end{align}
where $\delta_{k}^{l}$ and $\delta_{k}^{r}$ are nonlinear weights. In the design of nonlinear weights, the resolution needs to be preserved as high as possible. In smooth regions, the optimal sixth-order accuracy reconstruction is given at interface $x_{i+1/2}$. A combination of cubic polynomials is provided when discontinuity is not at interface $x_{i+1/2}$, and a combination of quadratic polynomials is given when discontinuity is just at interface $x_{i+1/2}$. Similar with the fifth-order WENO-JS and WENO-Z scheme \cite{WENO-JS}, the non-linear weights for the current WENO scheme are defined as
\begin{align*}
\delta^{l,r,JS}_{k}=\frac{\alpha^{l,r,JS}_{k}}{\sum_{m=1}^{m=4}\alpha^{l,r,JS}_{m}},~~\alpha^{l,r,JS}_{k}=\frac{d^{l,r}_{k}}{(\beta^{l,r}_{k}+\epsilon)^{p}},
\end{align*}
and
\begin{align*}
\delta^{l,r,Z}_{k}=&\frac{\alpha^{l,r,Z}_{k}}{\sum_{m=1}^{m=4}\alpha^{l,r,Z}_{m}}, ~~\alpha^{l,r,Z}_{k}=d^{l,r}_k\Big[1+\big(\frac{\tau}{\beta_k^{l,r}+\epsilon}\big)^2\Big],
\end{align*}
where $\tau$ is the global higher order reference value, which will be given in the following section. $\beta^{l,r}_{k}$ is the smooth indicator and calculated as the classical definition in \cite{WENO-JS}
\begin{align*}
\beta_k^{l,r}=\sum_{p=1}^{p_k-1}\Delta x^{2p-1}\int_{x_{i-1/2}}^{x_i+1/2}\big(\frac{\text{d}^p}{\text{d}x^p}w^{l,r}_k(x)\big)^2dx,
\end{align*}
where the smooth indicator $\beta^{l,r}_{1}$ corresponding to $S_{1,i+1/2}^{l,r}$ is replaced by that of the cubic polynomial $\widetilde{w}^{l,r}_1(x)$ on the stencil $\widetilde{S}^{l}_{1,i+1/2}=\{I_{i-2},I_{i-1},I_{i},I_{i+1}\}$ and $\widetilde{S}^{r}_{1,i+1/2}=\{I_{i},I_{i+1},I_{i+2},I_{i+3}\}$, and $p_0=2, p_1=p_2=p_3=3$. The following properties are satisfied for the non-linear weights
\begin{align*}
\begin{cases}
\text{discontinuity at~} x_{i-3/2},~~~\delta^{l,r}_{k}\ll O(1),~~k=0,1,\\
\text{discontinuity at~} x_{i-1/2},~~~\delta^{l,r}_{k}\ll O(1),~~k=0,1,2,\\
\text{discontinuity at~} x_{i+1/2},~~~\delta^{l,r}_{k}\ll O(1),~~k=1,2,3,\\
\text{discontinuity at~} x_{i+3/2},~~~\delta^{l,r}_{k}\ll O(1),~~k=2,3,\\
\text{discontinuity at~} x_{i+5/2},~~~\delta^{l,r}_{k}\ll O(1),~~k=3.
\end{cases}
\end{align*}
The details of smooth indicator for the left side of interface can be written as
\begin{align*}
\beta_0^l&=\frac{1}{4}(W_{i-2}-4W_{i-1}+3W_{i})^{2}+\frac{13}{12}(W_{i-2}-2W_{i-1}+W_{i})^{2},\\
\beta_1^l&=\frac{1}{36}(W_{i-2}-6W_{i-1}+3W_{i}+2W_{i+1})^{2}+ \frac{13}{12}(W_{i-1}-2W_{i}+W_{i+1})^{2} \\
&+\frac{1043}{960}(-W_{i-2}+3W_{i-1}-3W_{i}+W_{i+1})^{2} \\
&+\frac{1}{432}(W_{i-2}-6W_{i-1}+3W_{i}+2W_{i+1})(-W_{i-2}+3W_{i-1}-3W_{i}+W_{i+1}),\\
\beta_2^l&=\frac{1}{36}(-2W_{i-1}-3W_{i}+6W_{i+1}-W_{i+2})^{2}+ \frac{13}{12}(W_{i-1}-2W_{i}+W_{i+1})^{2} \\
&+\frac{1043}{960}(-W_{i-1}+3W_{i}-3W_{i+1}+W_{i+2})^{2} \\
&+\frac{1}{432}(-2W_{i-1}-3W_{i}+6W_{i+1}-W_{i+2})(-W_{i-1}+3W_{i}-3W_{i+1}+W_{i+2}),\\
\beta_3^l&=\frac{1}{36}(-11W_{i}+18W_{i+1}-9W_{i+2}+2W_{i+3})^{2}+ \frac{13}{12}(2W_{i}-5W_{i+1}+4W_{i+2}+W_{i+3})^{2} \\
&+\frac{1043}{960}(-W_{i}+3W_{i+1}-3W_{i+2}+W_{i+3})^{2} \\
&+\frac{1}{432}(-11W_{i}+18W_{i+1}-9W_{i+2}+2W_{i+3})(-W_{i}+3W_{i+1}-3W_{i+2}+W_{i+3}).
\end{align*}
According to the symmetry property, $\beta^{r}_{k}$ corresponding to the right side reconstruction can be obtained as well.

\subsection{Accuracy of the nonlinear schemes}
In this section, the accuracy of non-linear new scheme is analysed.  In the smooth region, the approximation error for the linear reconstruction can be written as
\begin{align*}
W(x_{i+1/2})-W^{opt}_{i+1/2}=A\Delta x^{6}+O(\Delta x^{7}),
\end{align*}
where $W(x_{i+1/2})$ is the exact solution at the interface $x_{i+1/2}$, $A$ is the Taylor expansion coefficient. With the sixth-order recursive WENO scheme, the reconstructed variables with nonlinear weights can be rewritten as
\begin{align}\label{order}
W^{l,r}_{i+1/2}=\sum_{k=0}^{k=3}d^{l,r}_{k}w^{l,r}_{k,i+1/2}+ \sum_{k=0}^{k=3}(\delta^{l,r}_{k}-d^{l,r}_{k})w^{l,r}_{k,i+1/2},
\end{align}
where $w^{l,r}_{k}(x)$ are the quadratic and cubic polynomials, and they approximate $W(x_{i+1/2})$ at least to $O(\Delta x^{3})$
\begin{align*}
w^{l,r}_{k,i+1/2}=W^{l,r}(x_{i+1/2})+B_{k}\Delta x^{3}+O(\Delta x^{4}).
\end{align*}
Substituting $w^{l,r}_{k,i+1/2}$ into the Eq.\eqref{order} and taking $\displaystyle\sum_{k=0}^{3}\delta^{l,r}_{k}=\sum_{k=0}^{3}d^{l,r}_{k}=1$ into account,  we have
\begin{align*}
W^{l,r}_{i+1/2}&=\sum_{k=0}^{k=3}d^{l,r}_{k}w^{l,r}_{k,i+1/2}+ \sum_{k=0}^{k=3}(\delta^{l,r}_{k}-d^{l,r}_{k})\big(W^{l,r}(x_{i+1/2})+B_{k}\Delta x^{3}+O(\Delta x^{4})\big) \\
&=W^{opt}_{i+1/2}+\Delta x^{3}\sum_{k=0}^{k=3}B_{k}(\delta^{l,r}_{k}-d^{l,r}_{k}) +\sum_{k=0}^{k=3}(\delta^{l,r}_{k}-d^{l,r}_{k})O(\Delta x^{4}),
\end{align*}
where the second and the third terms are the nonlinear remainders. In order to achieve the sixth-order of accuracy for the spatial discretization, the following equation condition needs to be satisfied \cite{Harten, WENO-M}
\begin{align*}
W(x_{i+1/2})-W^{l,r}_{i+1/2}=A^{'}\Delta x^{6}+O(\Delta x^{7}),
\end{align*}
where $A^{'}$ is a bounded variable satisfying Lipschitz continuity. Thus, the following sufficient condition is proposed for the nonlinear weights
\begin{align}\label{weight condition}
\delta^{l,r}_{k}-d^{l,r}_{k}=O(\Delta x^4).
\end{align}

For the WENO-JS weighting approach,  the sufficient condition \eqref{weight condition} can not be satisfied. Similar with the classical methodology, WENO-Z weighting approach is considered. Taylor expansion for the smooth indicators $\beta^{l}_{k}, k=1,2,3$ can be written as
\begin{align*}
\beta^{l}_{1}=(W_{i}^{(1)}\Delta x)^{2}+(\frac{13}{12}(W^{(2)}_{i})^{2}+\frac{1}{12}W_{i}^{(1)}W_{i}^{(3)})\Delta x)^{4}+\frac{1}{8}W_{i}^{(1)}W_{i}^{(4)}\Delta x^{5}+\Delta x^{6}, \\
\beta^{l}_{2}=(W_{i}^{(1)}\Delta x)^{2}+(\frac{13}{12}(W^{(2)}_{i})^{2}+\frac{1}{12}W_{i}^{(1)}W_{i}^{(3)})\Delta x)^{4}-\frac{1}{8}W_{i}^{(1)}W_{i}^{(4)}\Delta x^{5}+\Delta x^{6}, \\
\beta^{l}_{3}=(W_{i}^{(1)}\Delta x)^{2}+(\frac{13}{12}(W^{(2)}_{i})^{2}+\frac{1}{12}W_{i}^{(1)}W_{i}^{(3)})\Delta x)^{4}+\frac{5}{8}W_{i}^{(1)}W_{i}^{(4)}\Delta x^{5}+\Delta x^{6}.
\end{align*}
With the following local reference smooth indicator
\begin{align*}
\tau^{l,r}=\frac{-3\beta^{l,r}_{1}+2\beta^{l,r}_{2}+\beta^{l,r}_{3}}{6},
\end{align*}
we have
\begin{align*}
\tau^{l,r}=O(\Delta x^{6}).
\end{align*}
With the definition of the WENO-Z weighting approach, the sufficient condition Eq.\eqref{weight condition} is satisfied
\begin{align*}
\delta^{l,r}_{k}=d^{l,r}_{k}\Big[1+(\frac{\tau^{l,r}}{\beta^{l,r}_{k}+\epsilon})\Big]=d^{l,r}_{k}\big[1+O(\Delta x^{4})\big].
\end{align*}

\section{Numerical tests}
In this section, the numerical scheme will be presented to validate the current recursive WENO reconstruction. In the computation, two kinds of temporal discretization are considered.
The first one is the classical third-order TVD Runge-Kutta method \cite{TVD-rk}
\begin{align*}
W_{i}^{(1)}&=W_{i}^{n}+\mathcal{L}(W^{(0)}), \\
W_{i}^{(2)}&=\frac{3}{4}W_{i}^{n}+\frac{1}{4}W_{i}^{(1)}+\frac{1}{4}\Delta t\mathcal{L}(W^{(1)}), \\
W_{i}^{n+1}&=\frac{1}{3}W_{i}^{n}+\frac{2}{3}W_{i}^{(2)}+\frac{2}{3}\Delta t\mathcal{L}(W^{(2)}),
\end{align*}
with the recursive WENO reconstruction, the leading truncation error for the scheme is $O(\Delta x^6+\Delta t^3)$. With a fixed CFL number $\Delta t\sim\Delta x$, the order of accuracy will reduces to $3$. In order to keep the six-order accuracy, a small $\Delta t$ need to be used. Another choice is the two-stage fourth-order time-accurate discretization, which was developed for Lax-Wendroff flow solvers \cite{fourth-GRP,fourth-GKS}, which can be written as follows
\begin{align*}
W^*=&W^n+\frac{1}{2}\Delta t\mathcal{L}_{i}(W^n)+\frac{1}{8}\Delta t^2\frac{\partial}{\partial t}\mathcal{L}(W^n),\\
W^{n+1}=W^n&+\Delta t\mathcal{L}(w^n)+\frac{1}{6}\Delta t^2\big(\frac{\partial}{\partial t}\mathcal{L}(W^n)+2\frac{\partial}{\partial t}\mathcal{L}(W^*)\big),
\end{align*}
where the temporal derivative of $\mathcal{L}$ can be provided by the Lax-Wendroff flow solvers.

In this paper, the HLLC approximate Riemann solver \cite{Riemann-appro} in the TVD Runge-Kutta method and gas-kinetic flow solver \cite{GKS-Xu2} will be used two-stage discretization. For the two-dimensional computation, the fifth-order WENO reconstruction is used in the tangential direction. For each flux, the Gaussian quadratures are used in the tangential direction.

\subsection{Accuracy tests}
In this case, the advection of density perturbation is tested for the order of accuracy, and the initial condition is set as follows
\begin{align*}
\rho(x)=1+0.2\sin(\pi x),\ \  U(x)=1,\ \ \  p(x)=1, x\in[0,2].
\end{align*}
The periodic boundary condition is adopted, and the analytic solution is
\begin{align*}
\rho(x,t)=1+0.2\sin(\pi(x-t)),\ \ \  U(x,t)=1,\ \ \  p(x,t)=1.
\end{align*}
In the computation, a uniform mesh with $N$ points is used.

\begin{table}[!h]
\begin{center}
\def\temptablewidth{0.95\textwidth}
{\rule{\temptablewidth}{0.5pt}}
\begin{tabular*}{\temptablewidth}{@{\extracolsep{\fill}}c|cc|cc|cc}
mesh &  6th-linear& $L^1$ error ~  &6th-JS &  $L^1$ error ~& 6th-Z & $L^1$ error \\
\hline
10  &  1.3156E-4 & ~~     &  2.0032E-4 & ~~     & 1.3160E-4 &  ~~     \\
20  &  1.6514E-5 & 2.9940 &  1.7924E-5 & 3.4823 & 1.6514E-5 &  2.9944 \\
40  &  2.0664E-6 & 2.9985 &  2.0841E-6 & 3.1044 & 2.0664E-6 &  2.9985 \\
80  &  2.5836E-7 & 2.9997 &  2.5853E-7 & 3.0110 & 2.5836E-7 &  2.9997 \\
160 &  3.2298E-8 & 2.9999 &  3.2299E-8 & 3.0008 & 3.2298E-8 &  2.9999
\end{tabular*}
{\rule{\temptablewidth}{0.5pt}}
\end{center}
\vspace{-4mm} \caption{\label{tab1} Advection of density perturbation: the $L^1$ error and order of accuracy for the sixth-order spatial reconstruction with linear weights, JS-nonlinear weights and Z-nonlinear weights with third-order Runge-Kutta temporal discretization with $\Delta t=0.2\Delta x$.}
\end{table}

\begin{table}[!h]
\begin{center}
\def\temptablewidth{0.95\textwidth}
{\rule{\temptablewidth}{0.5pt}}
\begin{tabular*}{\temptablewidth}{@{\extracolsep{\fill}}c|cc|cc|cc}
mesh &  6th-linear& $L^1$ error ~  &6th-JS &  $L^1$ error ~& 6th-Z & $L^1$ error \\
\hline
10  &  1.0840E-05 &  ~     &  7.7952E-05 &  ~      & 1.0845E-05 &   ~     \\
20  &  1.7360E-07 & 5.9645 &  1.4662E-06 &  5.7324 & 1.7360E-07 &  5.9651 \\
40  &  2.7245E-09 & 5.9936 &  1.8252E-08 &  6.3279 & 2.7245E-09 &  5.9936 \\
80  &  4.2946E-11 & 5.9873 &  1.7202E-10 &  6.7293 & 4.2940E-11 &  5.9875 \\
160 &  5.6930E-12 & 2.9153 &  5.6987E-12 &  4.9158 & 5.6988E-12 &  2.9136
\end{tabular*}
{\rule{\temptablewidth}{0.5pt}}
\end{center}
\vspace{-4mm} \caption{\label{tab2} Advection of density perturbation: the $L^1$ error and order of accuracy for the sixth-order spatial reconstruction with linear weights, JS-nonlinear weights and Z-nonlinear weights with third-order Runge-Kutta temporal discretization with $\Delta t=\Delta x^2$.}
\end{table}

The $L^1$ and $L^2$ errors and orders with sixth-order spatial reconstruction and third-order Runge-Kutta method is present in Tab.\ref{tab1} for $\Delta t=0.2\Delta x$, and in Tab.\ref{tab2} for $\Delta t=\Delta x^2$.  For the time step with $\Delta t=\Delta x^2$, much more time steps is needed and more computational errors will be accumulated when the mesh is refined to $N=160$. To improve the temporal accuracy, the gas-kinetic scheme with the two-stage fourth-order time-accurate discretization and sixth-order spatial reconstruction is also tested with a fixed CFL number $CFL=0.2$ for different meshes. The $L^1$ errors and orders at $t=2$ are for sixth-order scheme with linear weights, JS-nonlinear weights and Z-nonlinear weights are presented in Tab.\ref{tab1}. With the mesh refinement to $N=640$, the expected orders of accuracy for different schemes.

\begin{table}[!h]
\begin{center}
\def\temptablewidth{0.95\textwidth}
{\rule{\temptablewidth}{0.5pt}}
\begin{tabular*}{\temptablewidth}{@{\extracolsep{\fill}}c|cc|cc|cc}
mesh &  6th-linear& $L^1$ error ~  &6th-JS &  $L^1$ error ~& 6th-Z & $L^1$ error \\
\hline
40  &   1.7087E-07 &   ~    &  5.9048E-07 &      ~      &   1.7144E-07 &   ~    \\
80  &   2.6805E-09 & 5.9942 &  1.3025E-08 &     5.5025  &   2.6826E-09 &  5.9979\\
160 &   4.1980E-11 & 5.9966 &  2.3328E-10 &     5.8030  &   4.1988E-11 &  5.9974\\
320 &   6.5954E-13 & 5.9921 &  3.8936E-12 &     5.9047  &   6.5877E-13 &  5.9940\\
640 &   1.6272E-14 & 5.3409 &  6.3352E-14 &     5.9415  &   1.7459E-14 &  5.2376
\end{tabular*}
{\rule{\temptablewidth}{0.5pt}}
\end{center}
\vspace{-4mm} \caption{\label{tab3} Advection of density perturbation: the $L^1$ error and order of accuracy for the gas-kinetic scheme for the sixth-order spatial reconstruction with linear weights, JS-nonlinear weights and Z-nonlinear weights.}
\end{table}

\begin{figure}[!htb]
\centering
\includegraphics[width=0.495\textwidth]{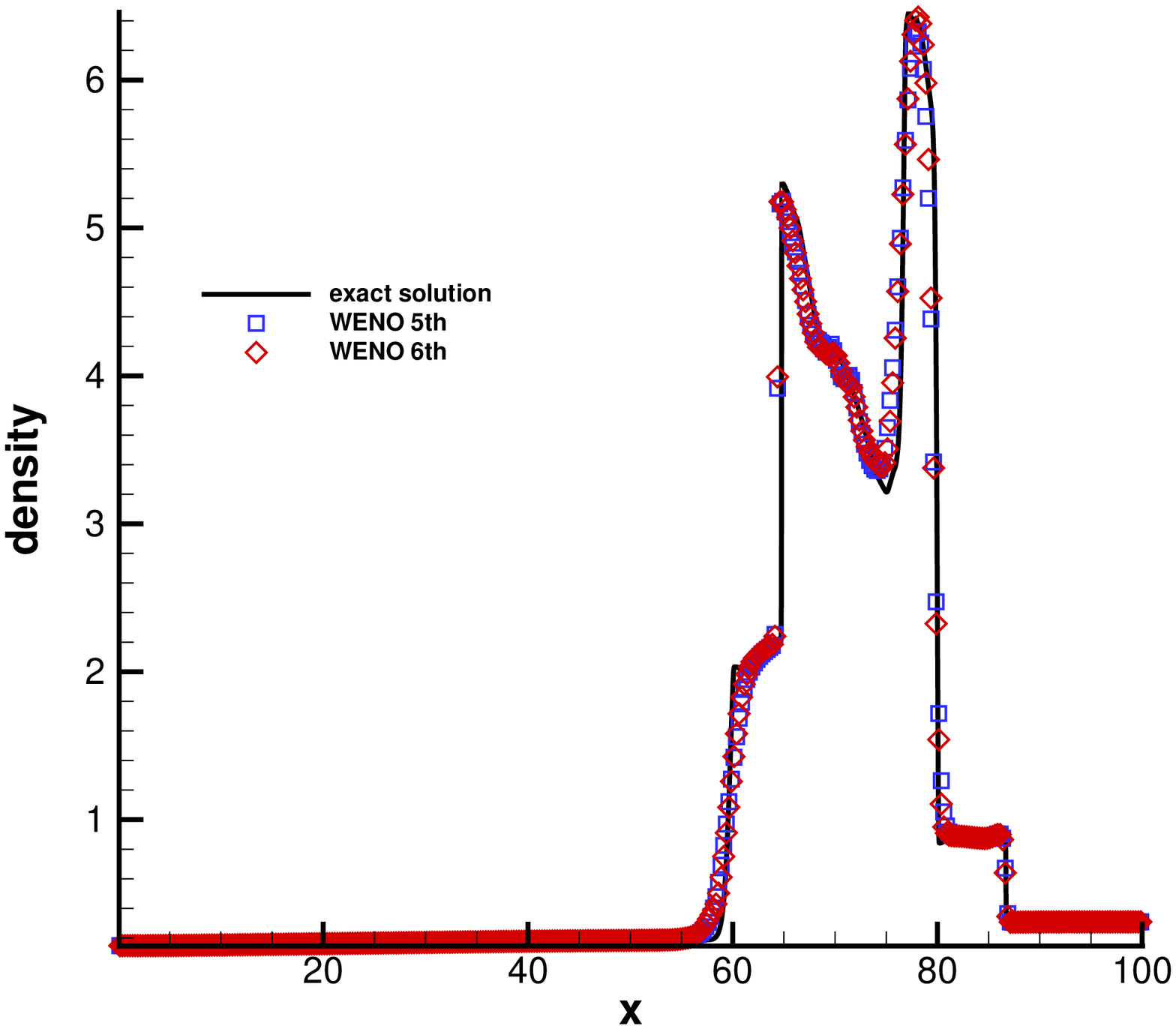}
\includegraphics[width=0.495\textwidth]{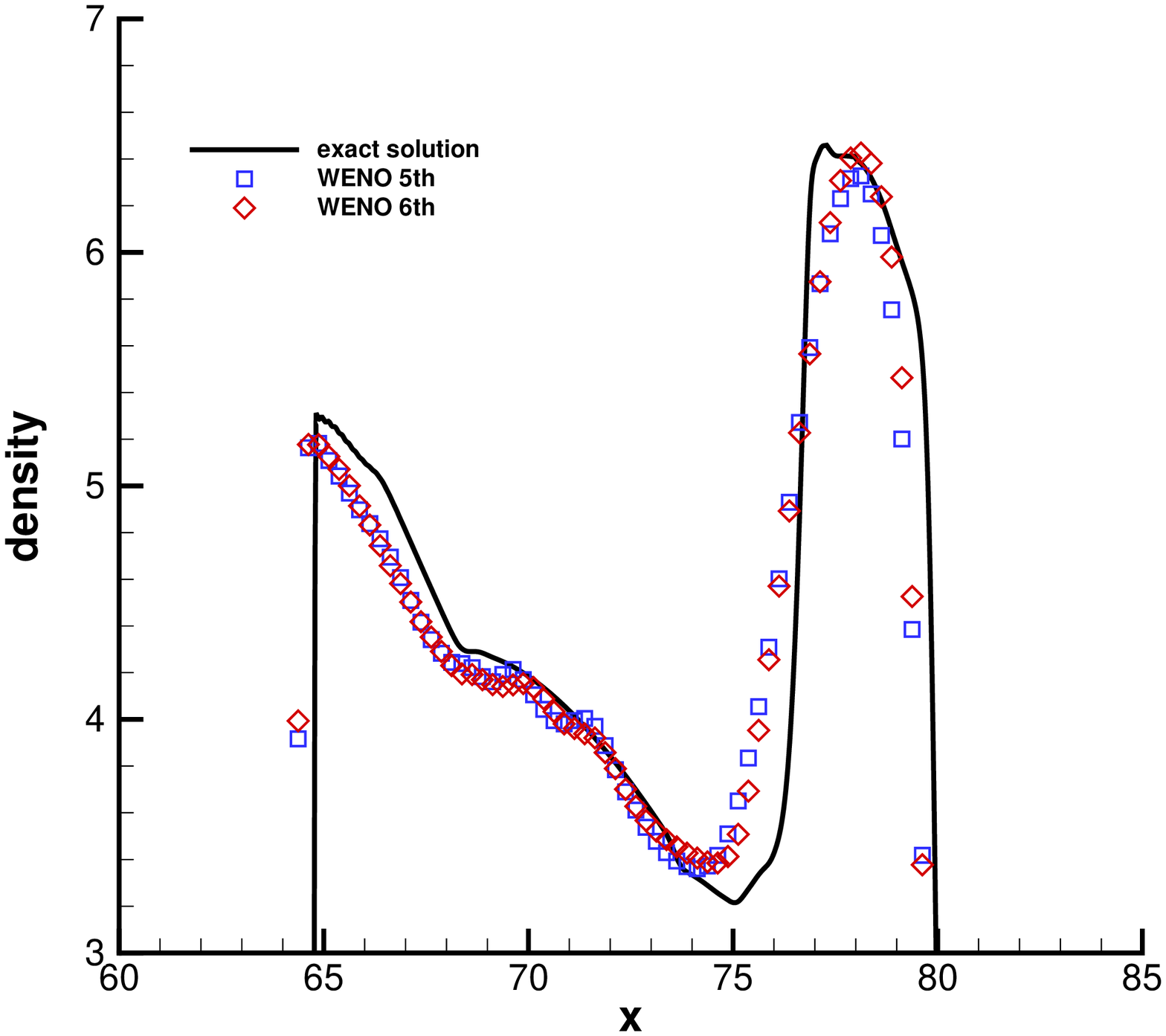}
\caption{\label{1d-riemann-1}1D Riemann problem: the density, velocity and pressure distributions for the blast-wave problem at $t=3.8$ with $400$ cells with fifth-order, sixth-order WENO-JS scheme.}
\end{figure}

\subsection{One dimensional Riemann problems}
For one-dimensional case, three Riemann problems are considered. The first one is the Woodward-Colella blast wave problem \cite{Case-Woodward}. The computational domain is $[0,100]$ with $400$ uniform mesh points. The reflected boundary conditions are imposed on both ends and the initial conditions are given as follows
\begin{equation*}
(\rho,U,p) =\left\{\begin{array}{ll}
(1, 0, 1000), \ \ \ \ & 0\leq x<10,\\
(1, 0, 0.01), & 10\leq x<90,\\
(1, 0, 100), &  90\leq x\leq 100.
\end{array} \right.
\end{equation*}
The density, velocity, and pressure distributions for the fifth-order, sixth-order WENO-JS scheme with and the exact solutions are presented in Fig.\ref{1d-riemann-1} for the blast wave problem at $t=3.8$. The numerical results agree well with the exact solutions. Especially, the sixth-order scheme resolves the local extreme values of blast-wave profile better than the fifth-order scheme.

\begin{figure}[!htb]
\centering
\includegraphics[width=0.495\textwidth]{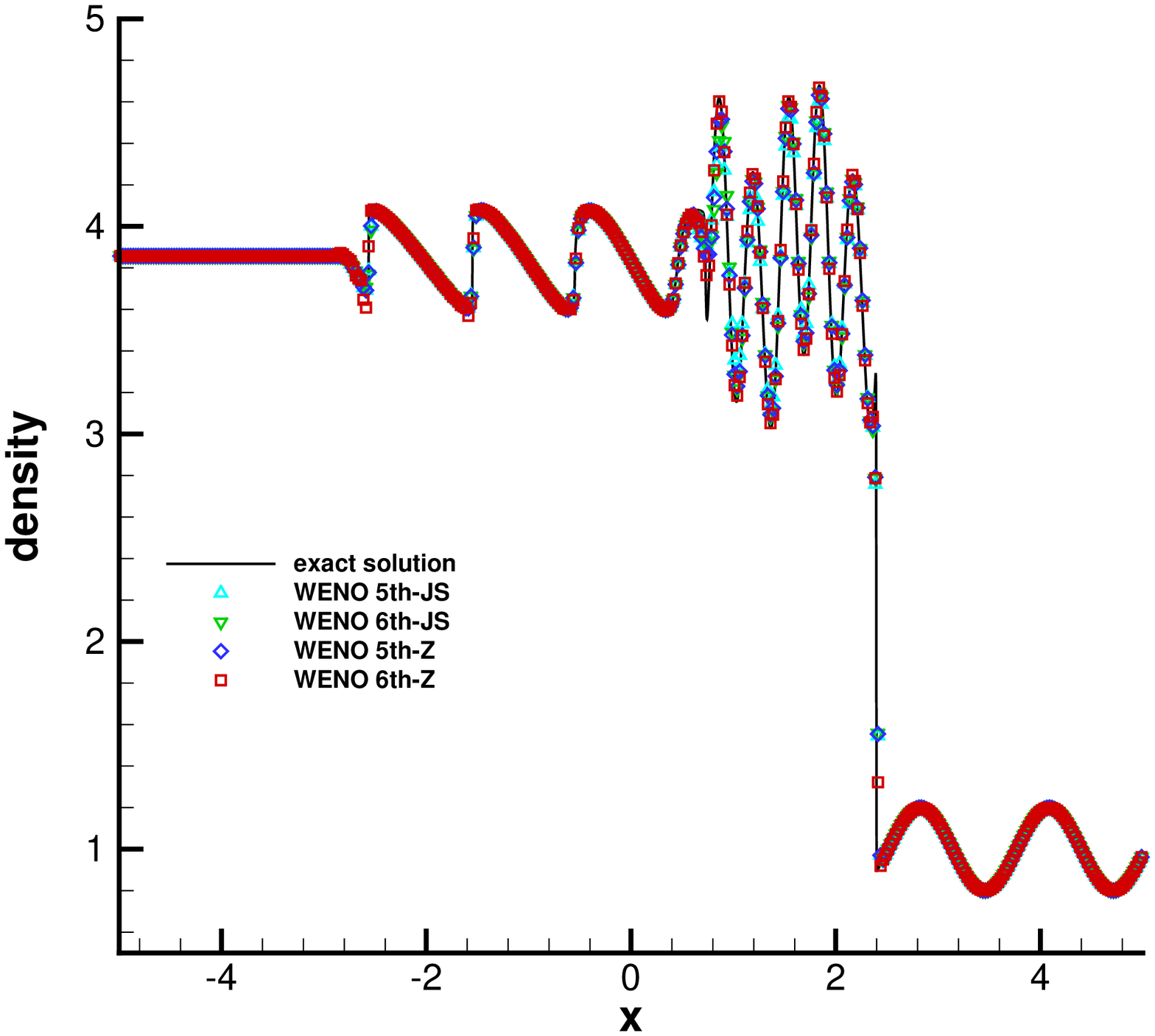}
\includegraphics[width=0.495\textwidth]{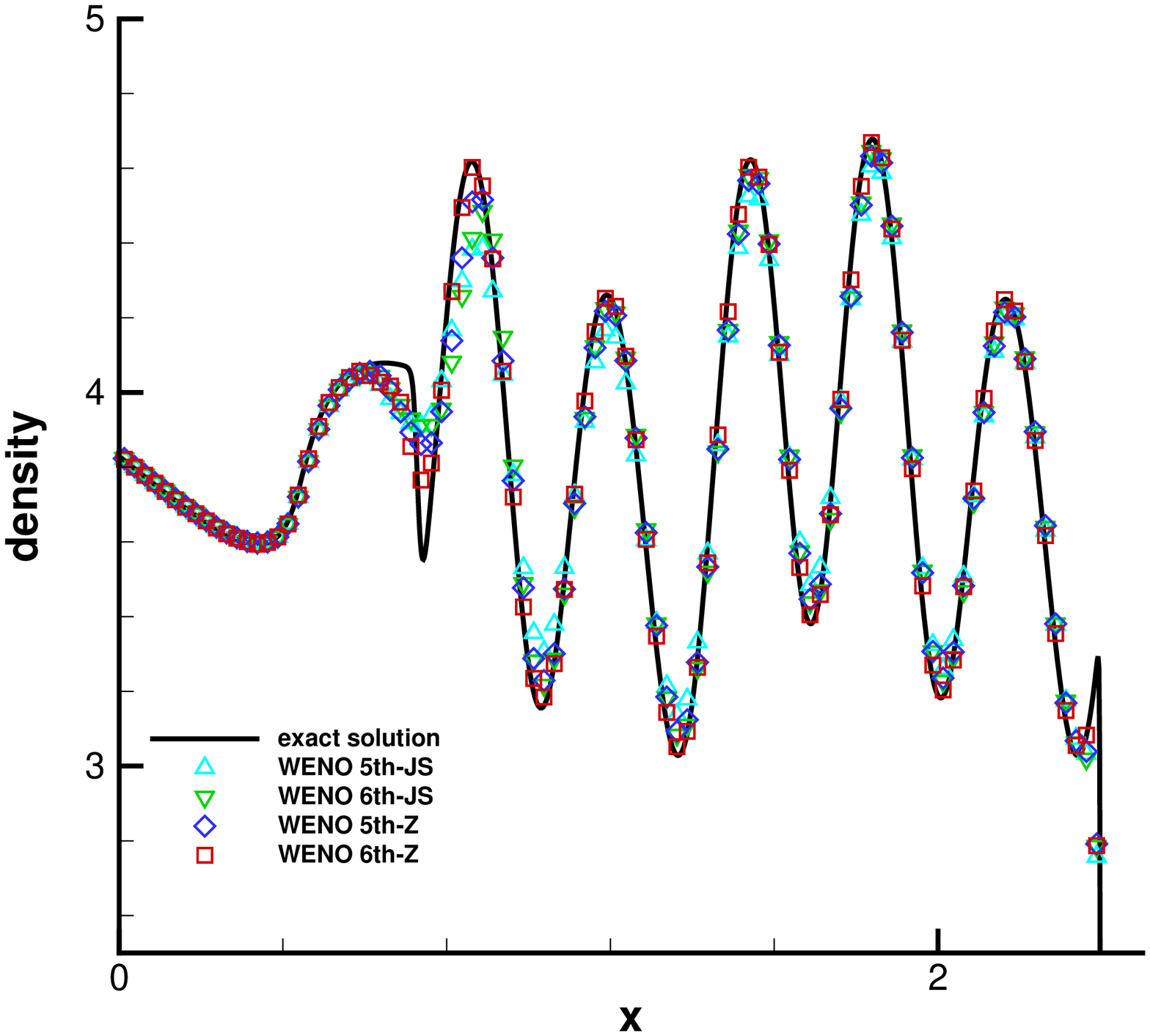}
\caption{\label{1d-riemann-2} 1D Riemann problem: the density distribution and local enlargement of Shu-Osher problem with fifth-order, sixth-order WENO-JS and WENO-Z scheme at $t=1.8$ with $1000$ cells.}
\end{figure}

\begin{figure}[!htb]
\centering
\includegraphics[width=0.495\textwidth]{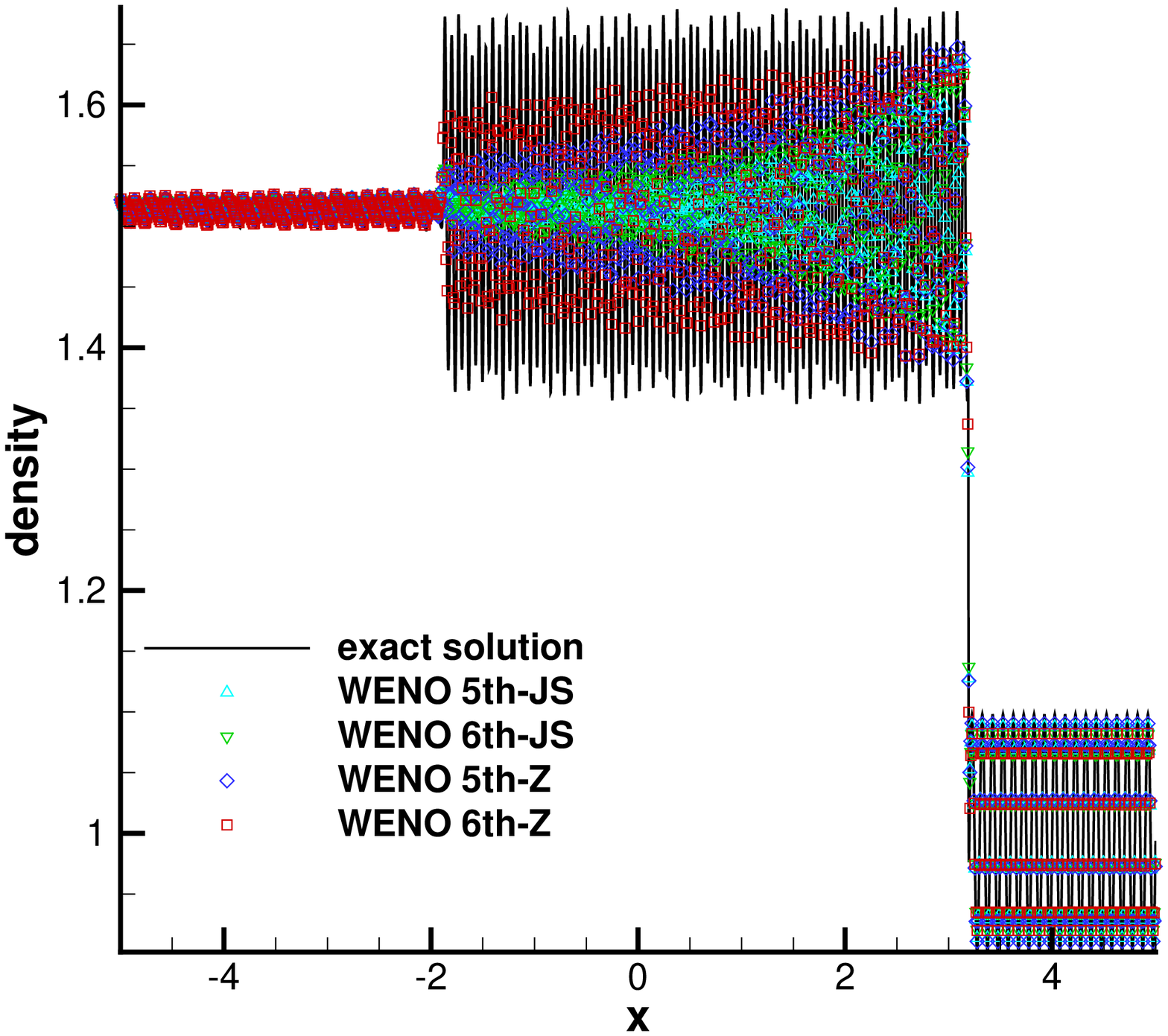}
\includegraphics[width=0.495\textwidth]{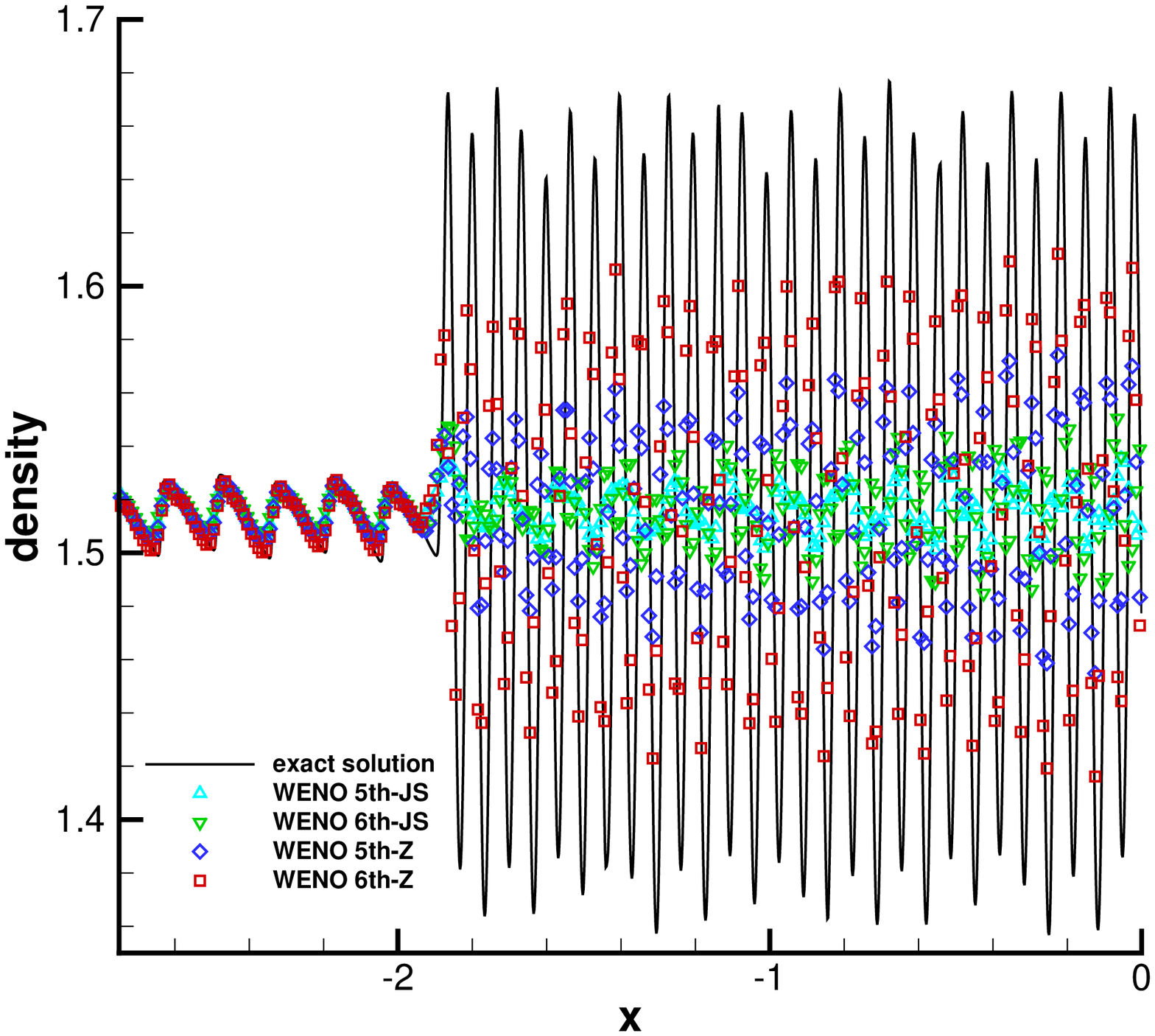}
\caption{\label{1d-riemann-3} 1D Riemann problem: the density distribution and local enlargement of Titarev-Toro problem with fifth-order, sixth-order WENO-JS and WENO-Z scheme at $t=5$ with $400$ cells.}
\end{figure}

In the one-dimensional case, another standard test case is the Shu-Osher problem \cite{Case-Shu-Osher}. The aim of this case is to test the ability of high-order numerical scheme to capture the high frequency waves. The computational domain is $[-5,5]$ and the flow field is initialized as
\begin{equation*}
(\rho,U,p)=\left\{\begin{array}{ll}
(3.857134, 2.629369, 10.33333),  \ \ \ \ &  x \leq -4,\\
(1 + 0.2\sin (5x), 0, 1),  &  -4 <x.
\end{array} \right.
\end{equation*}
As the extension of the Shu-Osher problem, the Titarev-Toro problem \cite{Case-Titarev-Toro} is a more severe test case for the oscillatory wave interacting with shock. The initial condition for this cases is given as follows
\begin{equation*}
(\rho,U,p)=\left\{\begin{array}{ll}
(1.515695,0.523346,1.805),  \ \ \ \ & -5< x \leq -4.5,\\
(1 + 0.1\sin (20\pi x), 0, 1),  &  -4.5 <x <5.
\end{array} \right.
\end{equation*}
The computed density profile and local enlargement for the Shu-Osher problem with $400$ uniform mesh points at $t =1.8$ and for the Titarev-Toro problem with $1000$ uniform mesh points at $t =5$ are shown in Fig.\ref{1d-riemann-2} and Fig.\ref{1d-riemann-3}, respectively. Fifth-order, sixth-order WENO-JS and WENO-Z scheme are used to test the performance of different orders and different nonlinear weights in the WENO reconstruction.
As analyzed in \cite{WENO-Z}, the numerical results with WENO-Z reconstruction resolves the local extreme values better than that from WENO-JS reconstruction. Due to the high order of accuracy for reconstruction, the sixth-order scheme resolves performs better than the fifth-order scheme.

\begin{figure}[!htb]
\centering
\includegraphics[width=0.8\textwidth]{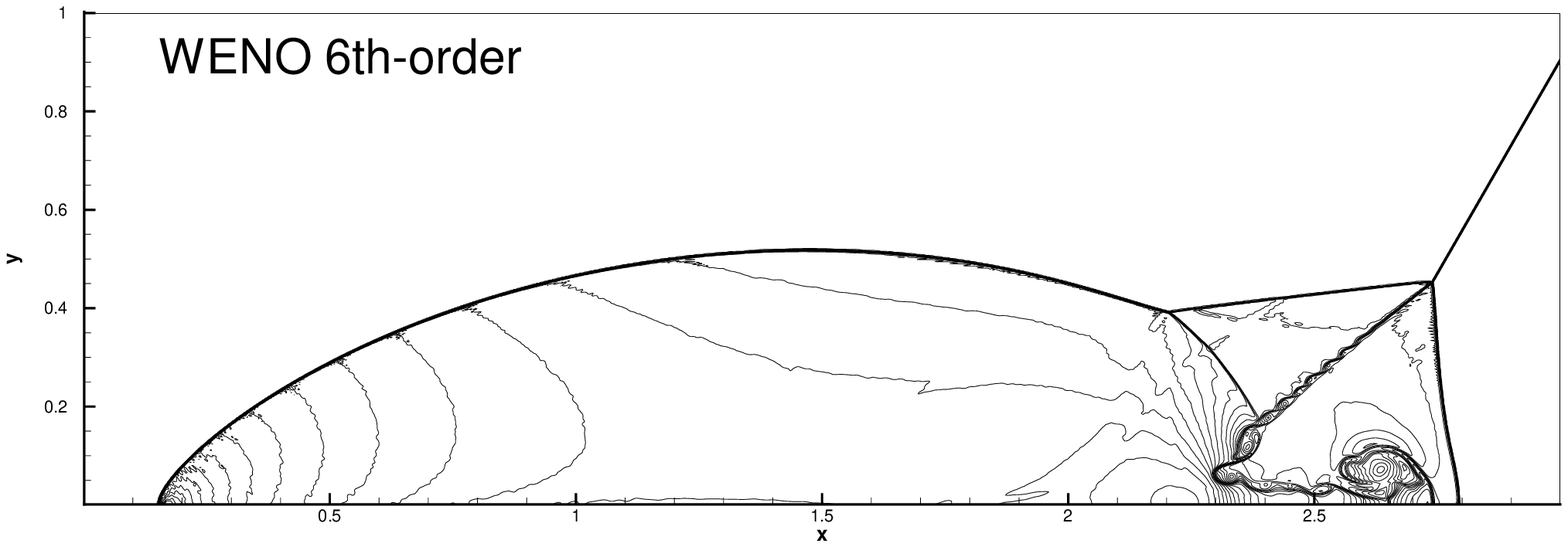}
\includegraphics[width=0.8\textwidth]{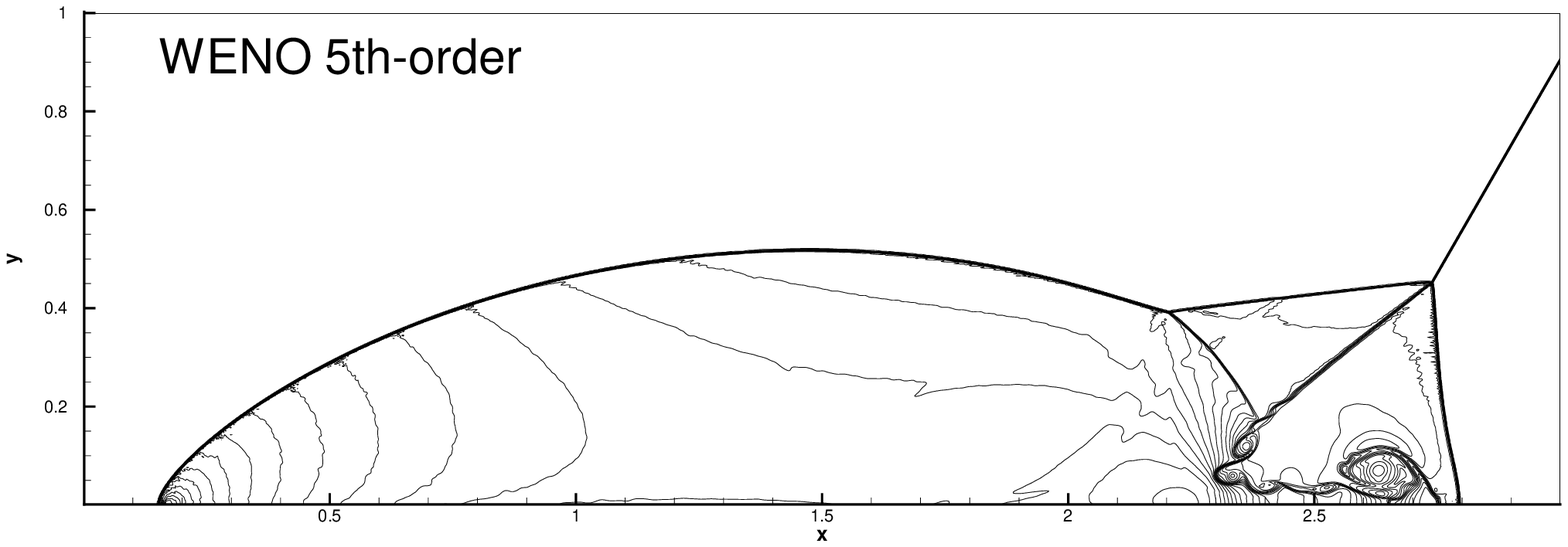}
\caption{\label{double-mach-1} Double Mach reflection: the density contours for sixth-order and fifth-order WENO reconstructions with  $1440\times480$ mesh points.}
\centering
\includegraphics[width=0.49\textwidth]{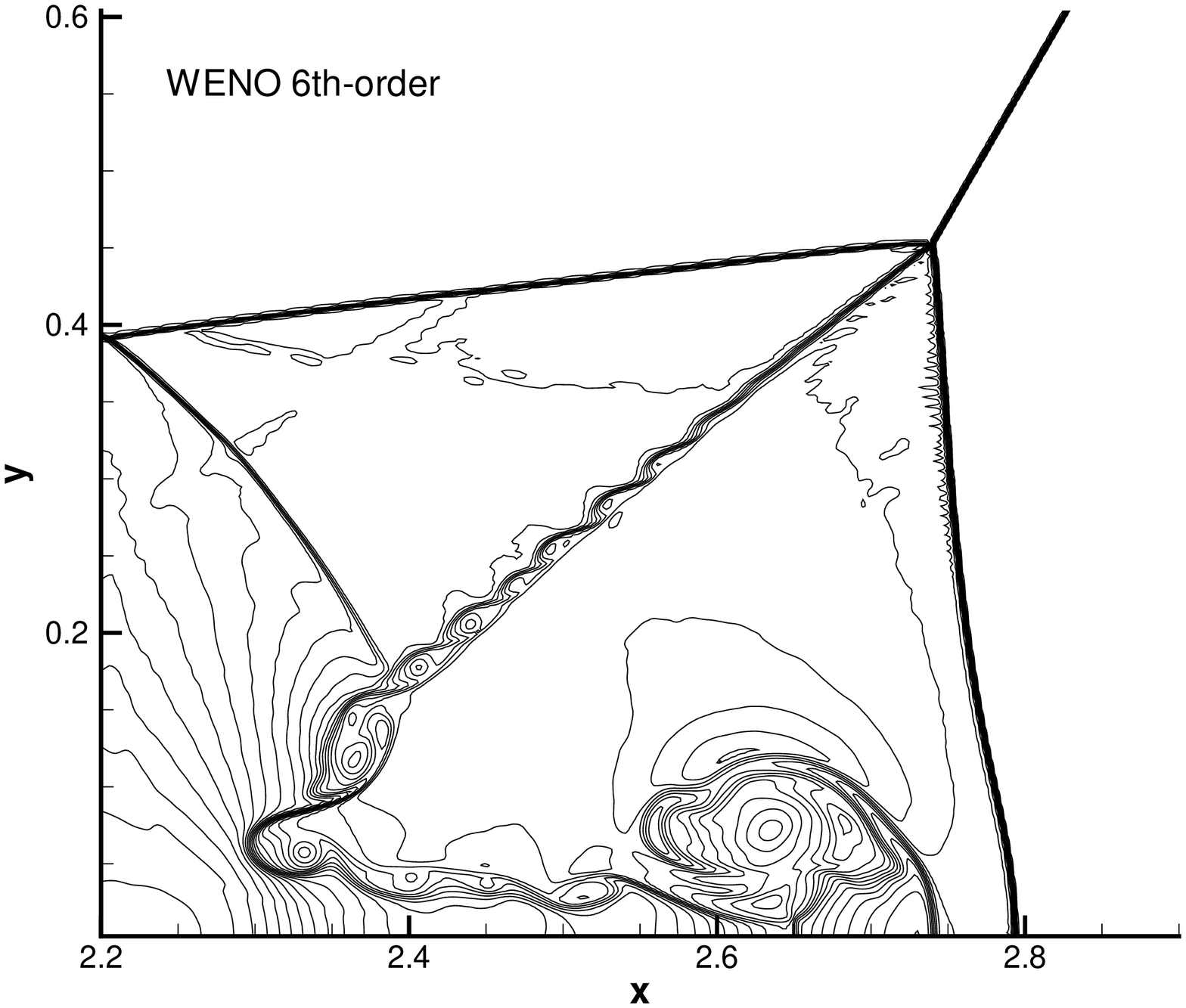}
\includegraphics[width=0.49\textwidth]{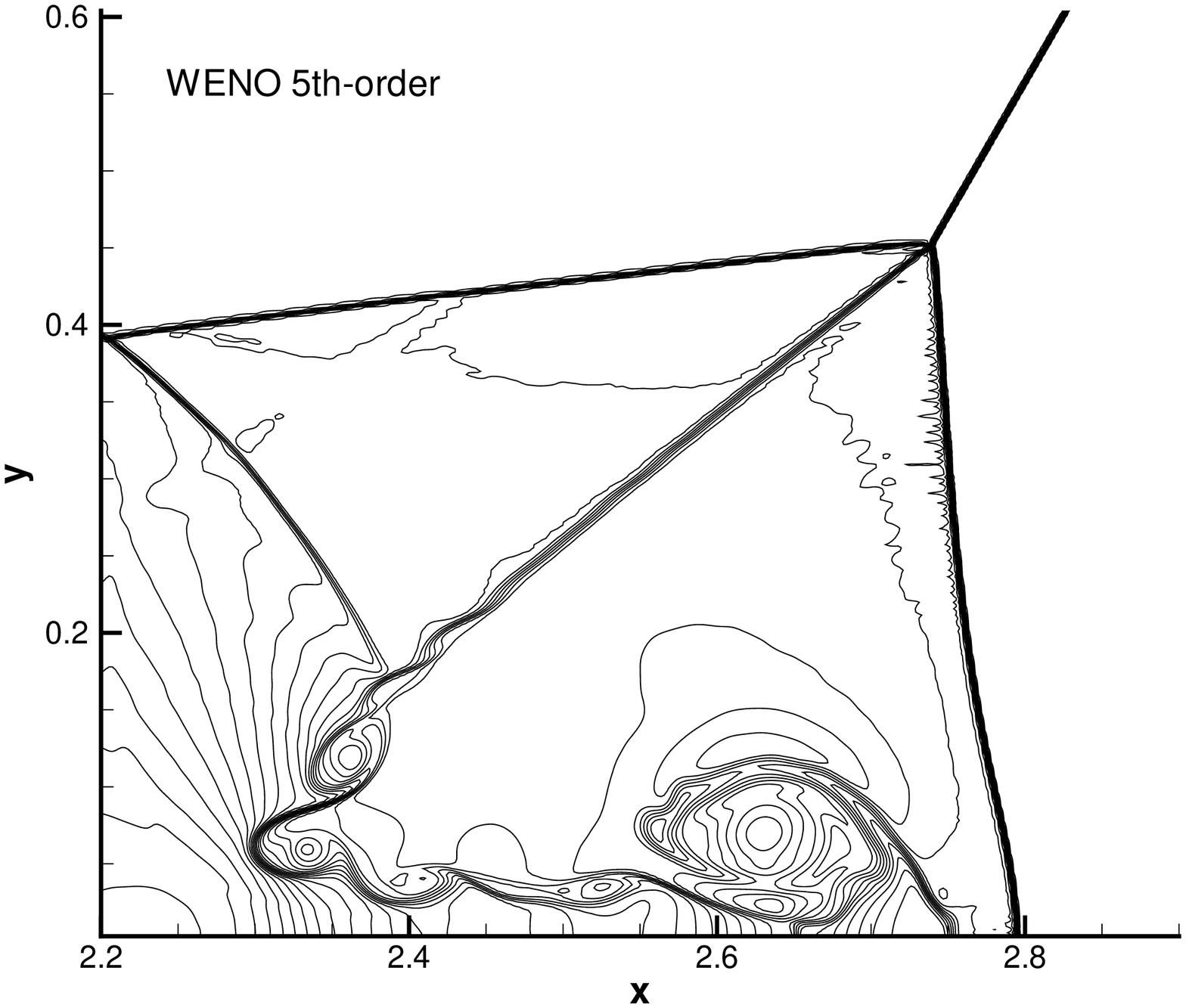}
\caption{\label{double-mach-2} Double Mach reflection: the local enlarged density distributions around the triple point for sixth-order and fifth-order WENO reconstructions with  $1440\times480$ mesh points.}
\end{figure}

\subsection{Double Mach reflection problem}
This problem was extensively studied by Woodward and Colella \cite{Case-Woodward} for the inviscid flow. The computational domain is $[0,4]\times[0,1]$, and a solid wall lies at the bottom of the computational domain starting from $x =1/6$. Initially a right-moving Mach 10 shock is positioned at $(x,y)=(1/6, 0)$, and makes a $60^\circ$ angle with the x-axis. The initial pre-shock and post-shock conditions are
\begin{align*}
(\rho, U, V, p)&=(8, 4.125\sqrt{3}, -4.125,
116.5),\\
(\rho, U, V, p)&=(1.4, 0, 0, 1).
\end{align*}
The reflective boundary condition is used at the wall, while for the rest of bottom boundary, the exact post-shock condition is imposed. At the top boundary, the flow variables are set to describe the exact motion of the Mach $10$ shock. Fifth-order and sixth-order scheme WENO-JS reconstructions are used in this case. The density distributions and local enlargement with $1440\times480$ uniform mesh points at $t=0.2$ for fifth-order and sixth-order scheme are shown in Fig.\ref{double-mach-1} and Fig.\ref{double-mach-2}, respectively. These two schemes resolve the flow structure under the triple Mach stem clearly. Compared with
the fifth-order scheme, the current sixth-order scheme is less dissipative and the instability of contact line are better resolved.

\begin{figure}[!htb]
\centering
\includegraphics[width=0.48\textwidth]{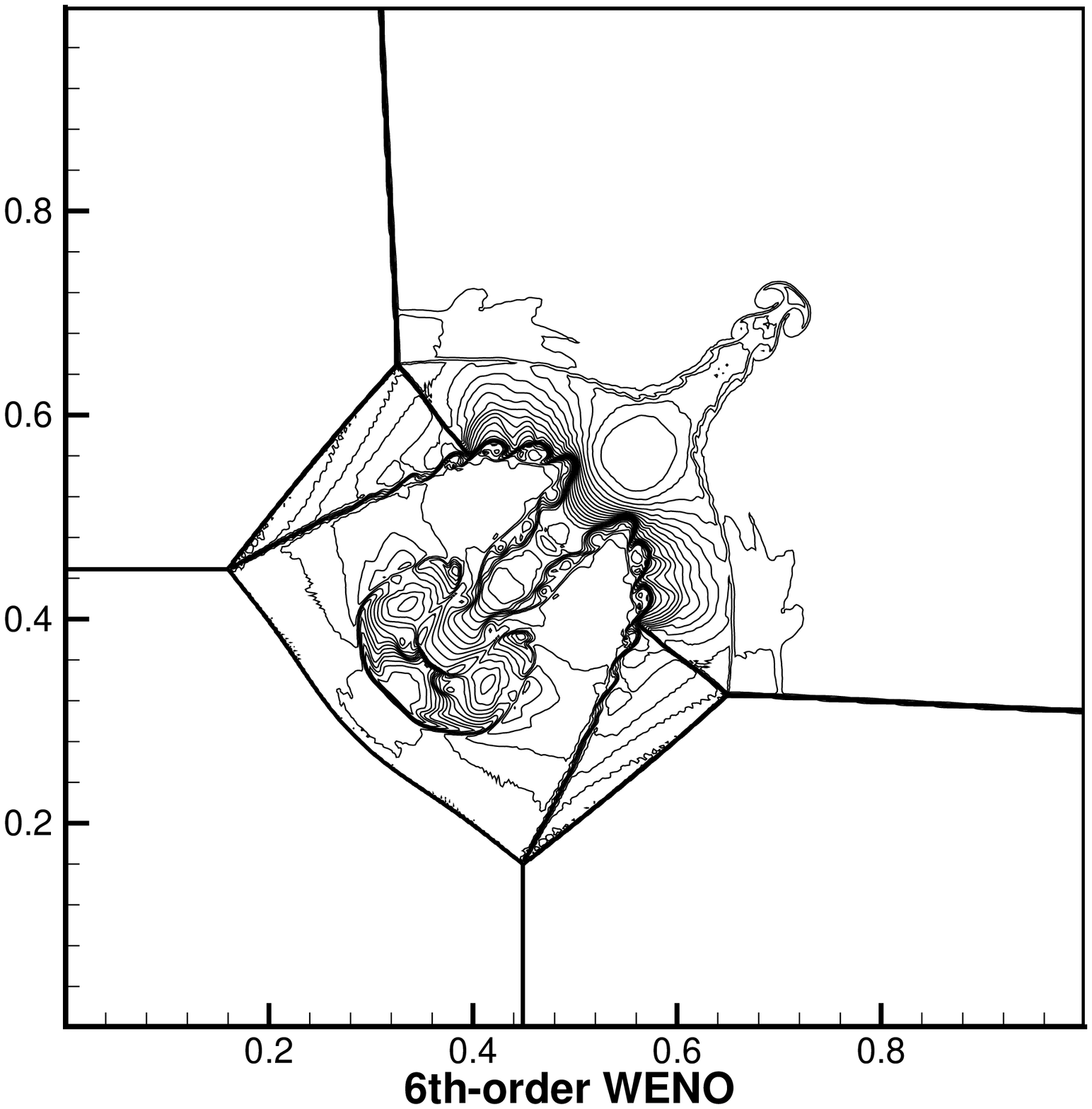}
\includegraphics[width=0.48\textwidth]{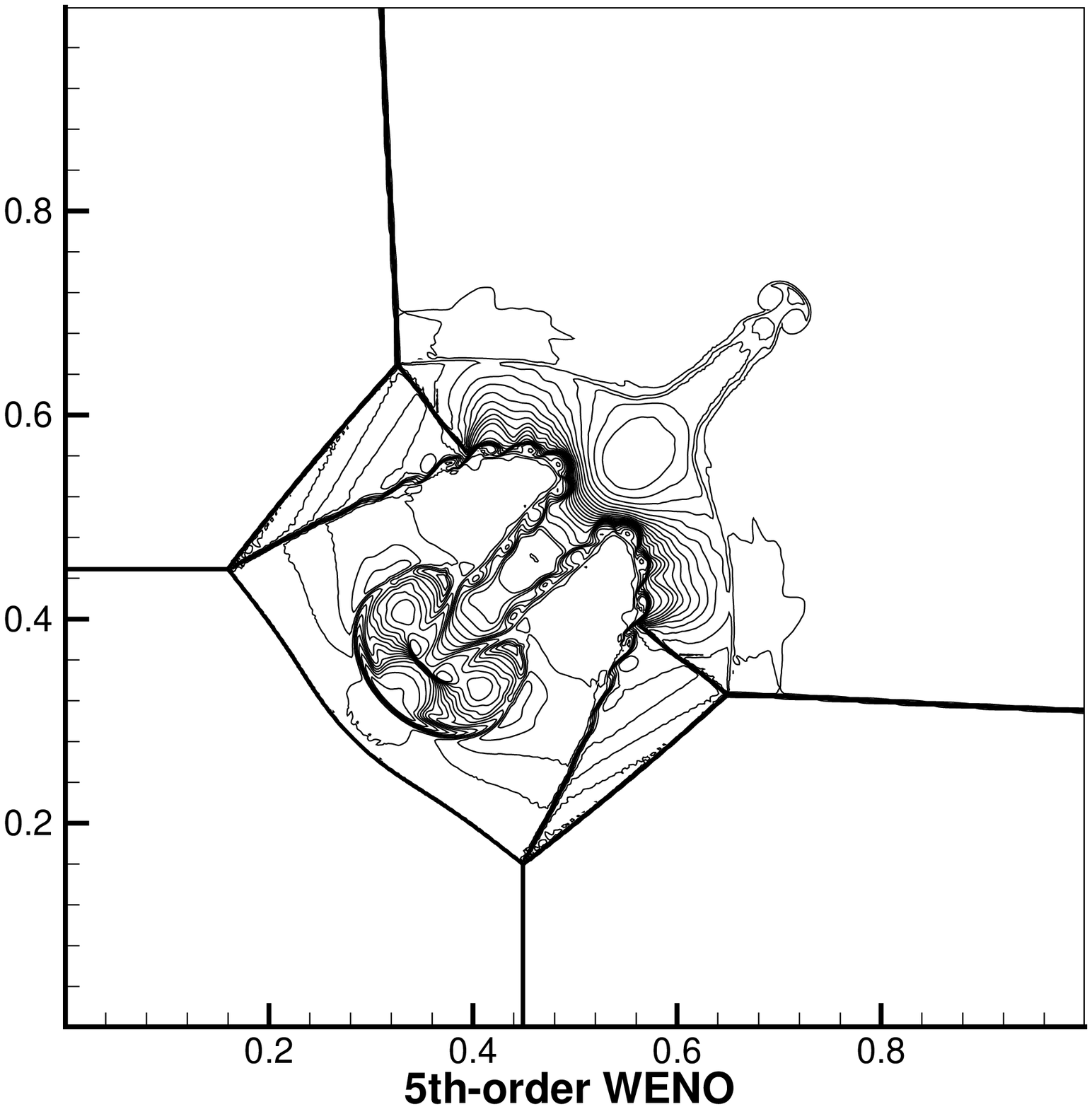}
\caption{\label{2d-riemann-1} 2D Riemann problem: the density distributions for sixth-order and fifth-order WENO reconstructions with  $500\times500$ mesh points for the interactions of shocks.}

\end{figure}

\begin{figure}[!htb]
\centering
\includegraphics[width=0.48\textwidth]{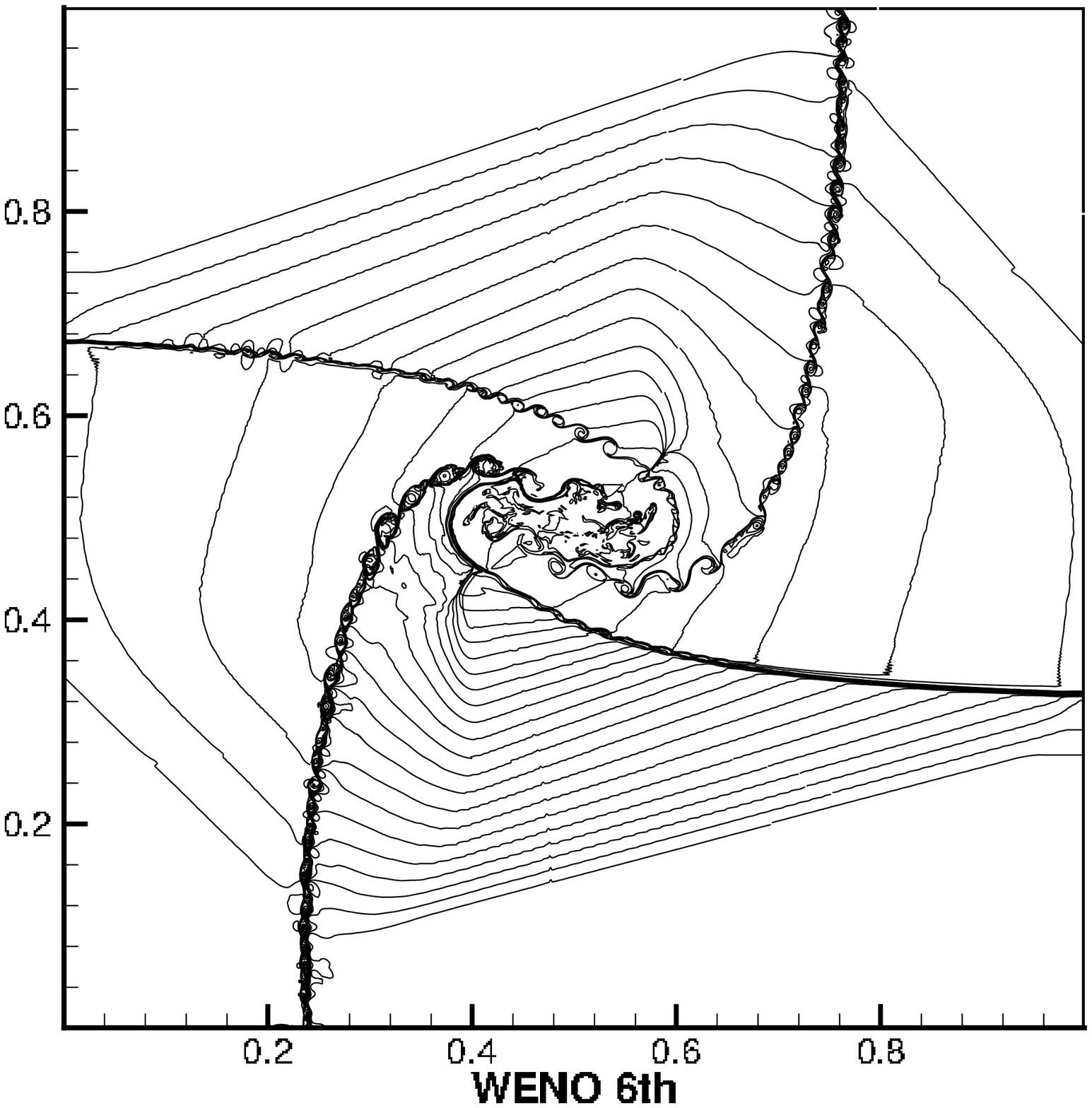}\includegraphics[width=0.48\textwidth]{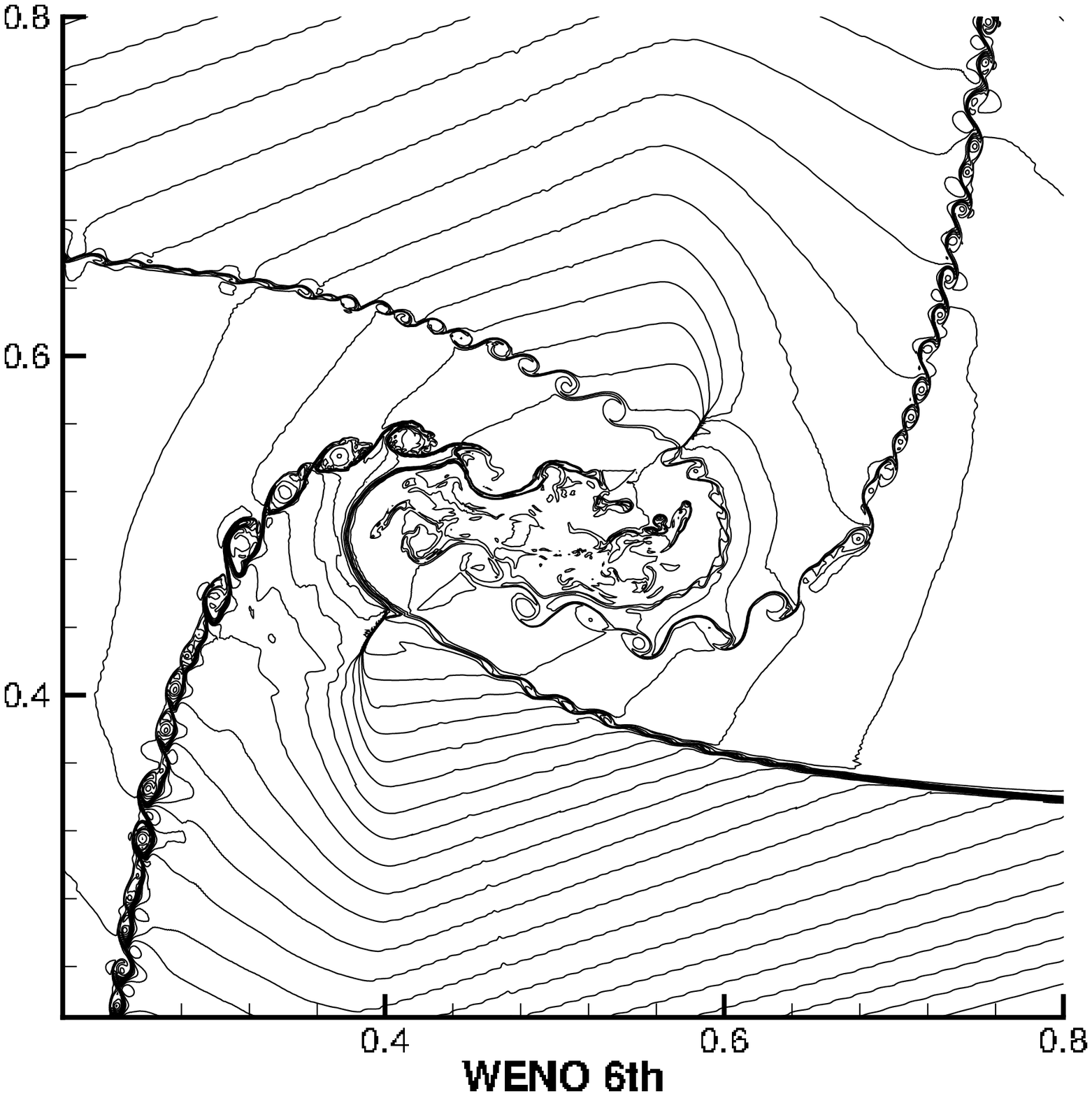}\\
\includegraphics[width=0.48\textwidth]{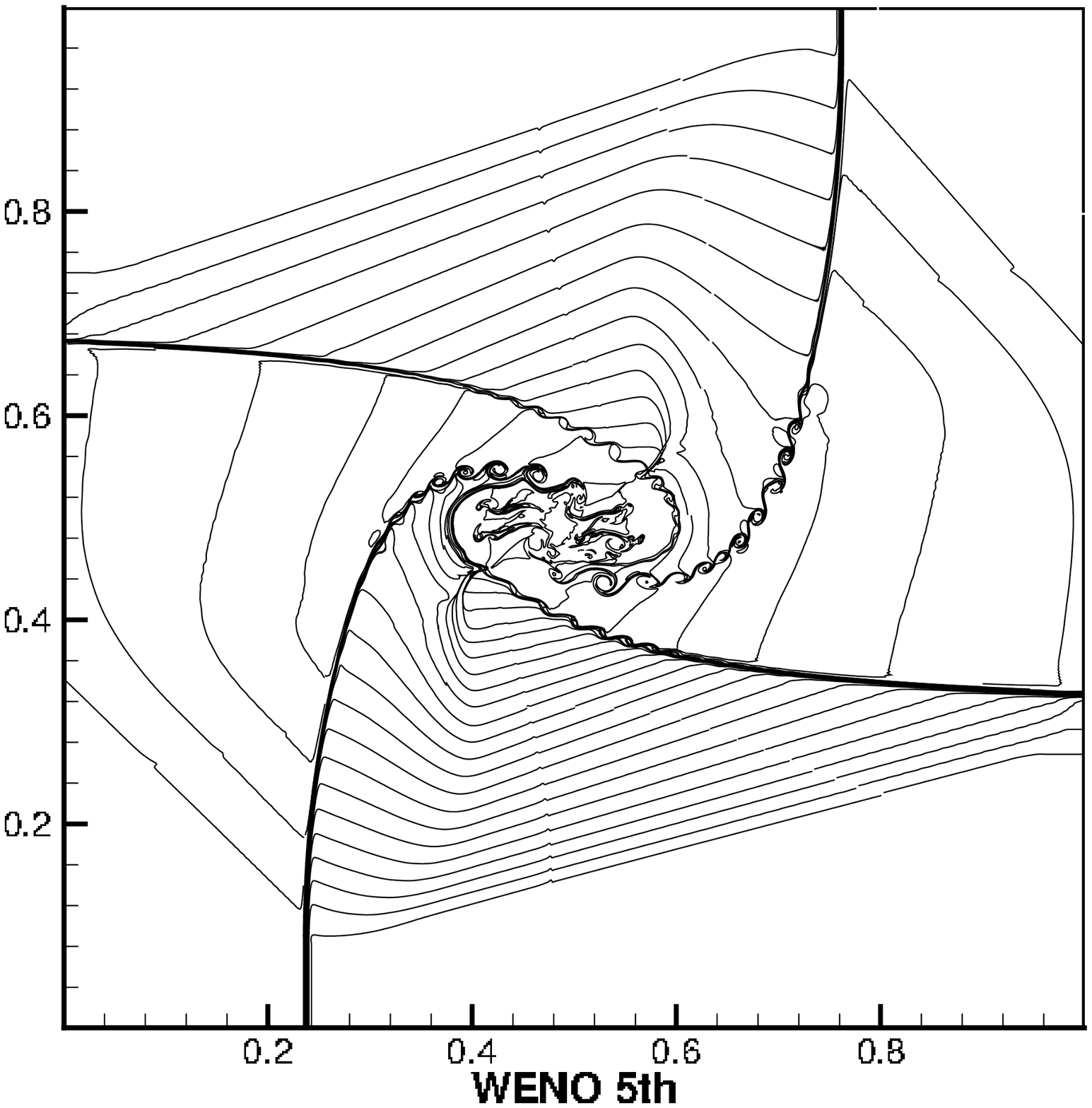}\includegraphics[width=0.48\textwidth]{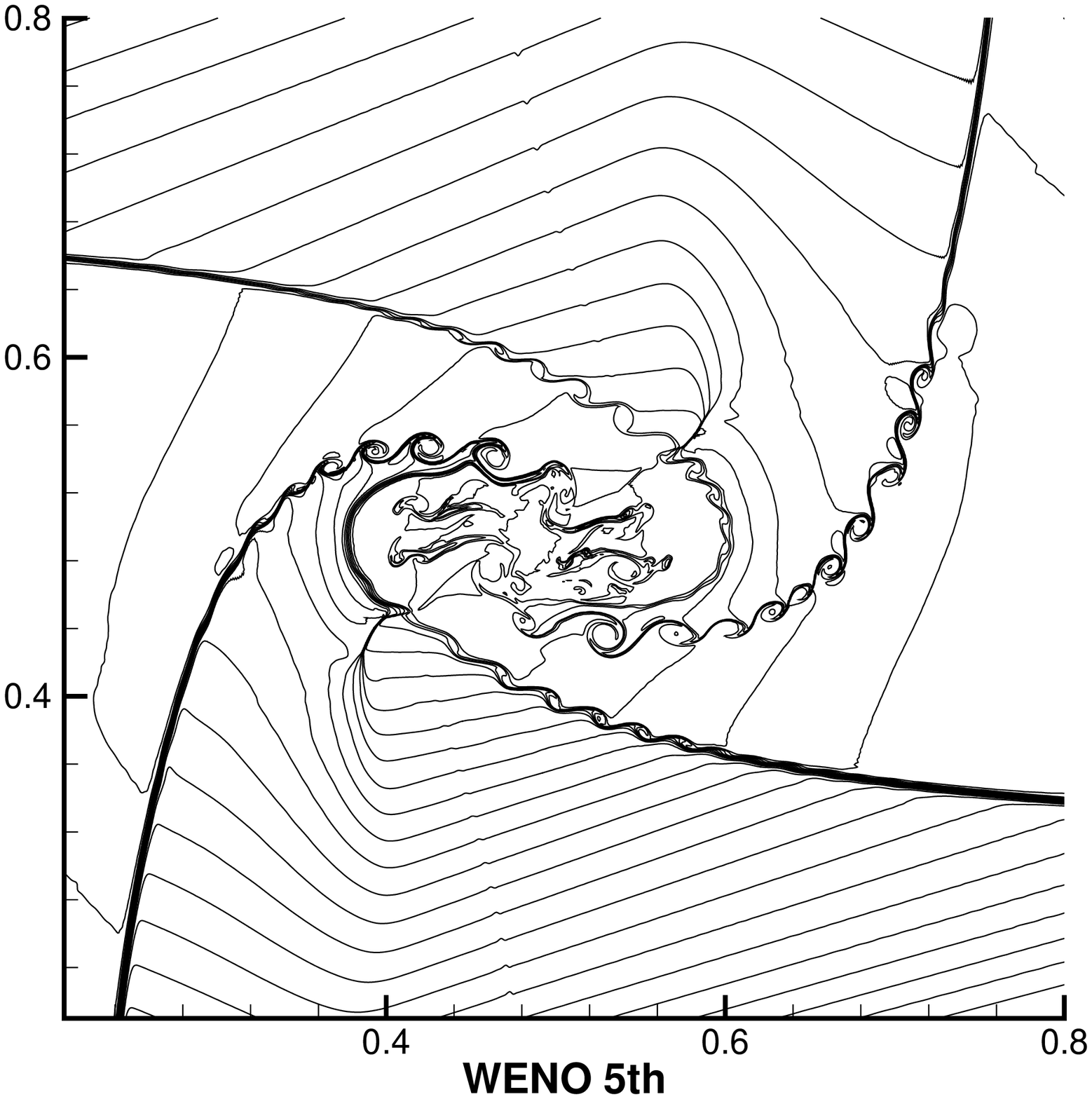}
\caption{\label{2d-riemann-2} 2D Riemann problem: the density distributions and local enlargement for sixth-order and fifth-order WENO reconstructions with  $1500\times1500$ mesh points for the interactions of contact discontinuities.}
\end{figure}

\subsection{Two-dimensional Riemann problems}
In this case, two examples of two-dimensional Riemann problems are considered, which involve the interactions of shocks, and the interaction of contact discontinuities \cite{Case-Riemann-1,Case-Riemann-2}.
In first case, the interaction of four shocks $\overleftarrow{S_{21}}\overleftarrow{S_{32}} \overleftarrow{S_{41}} \overleftarrow{S_{34}}$ is tested, where the backward rarefaction wave connecting the areas $\Omega_l$ and $\Omega_r$  are denoted as $\overleftarrow{S_{lr}}$. To obtain the detailed flow structure with less computational mesh points, the initial conditions are given as follows
\begin{equation*}
(\rho,U,V,p)=
\left\{\begin{array}{ll}
(1.5,0,0,1.5) \ \ \ &\Omega_1: x>0.7,y>0.7,\\
(0.5323,1.206,0,0.3), &\Omega_2: x<0.7,y>0.7,\\
(0.138,1.206,1.206,0.029), &\Omega_3: x<0.7,y<0.7,\\
(0.5323,0,1.206,0.3), &\Omega_4: x>0.7,y<0.7.
\end{array} \right.
\end{equation*}
The computational domain is $[0,1]\times[0,1]$, and the non-reflecting boundary conditions are used in all boundaries. The numerical solution is given in Fig.\ref{2d-riemann-1} at $t=0.6$, where the uniform mesh with $\Delta x=\Delta y=1/500$ is used. This case is just the mathematical formation of the double Mach problem \cite{Case-Woodward} and the symmetric line $x=y$ can be regarded as the rigid wall. The sixth-order and fifth-order schemes with WENO-JS reconstructions are tested to simulate the wave patters resulting from the interaction of shocks. The small scaled vortices are resolved  sharply using the current scheme, and the sixth-order scheme is less dissipative than the fifth-order one.

In second case, the interaction of four contact discontinuities $J_{21}^-J_{32}^- J_{41}^- J_{34}^-$ is tested, where the backward contact discontinuities connecting the areas $\Omega_l$ and $\Omega_r$  are denoted as $J_{lr}^-$.
The initial conditions for the this case are given as follows
\begin{equation*}
(\rho,U,V,p)=\left\{\begin{aligned}
         &(1 ,0.75,-0.5,1), &\Omega_1: x>0.5,y>0.5,\\
         &(2,0.75,0.5,1), &\Omega_2: x<0.5,y>0.5,\\
         &(1,-0.75,0.5,1), &\Omega_3: x<0.5,y<0.5,\\
         &(3,-0.75,-0.5,1), & \Omega_4: x>0.5,y<0.5.
                          \end{aligned} \right.
                          \end{equation*}
Their instantaneous interaction results in an entropy wave and a vortex sheet. The computational domain is $[0,1]\times[0,1]$, and the non-reflecting boundary conditions are also used in all boundaries.  The sixth-order and fifth-order schemes with WENO-JS reconstructions are tested. To obtain the detailed flow structure, the uniform mesh with $\Delta x=\Delta y=1/1500$ are used, and the numerical solution is given in Fig.\ref{2d-riemann-2} at $t=0.35$.  More small scaled vortices are resolved sharply by the sixth-order scheme.

\begin{figure}[!htb]
\centering
\includegraphics[width=0.25\textwidth]{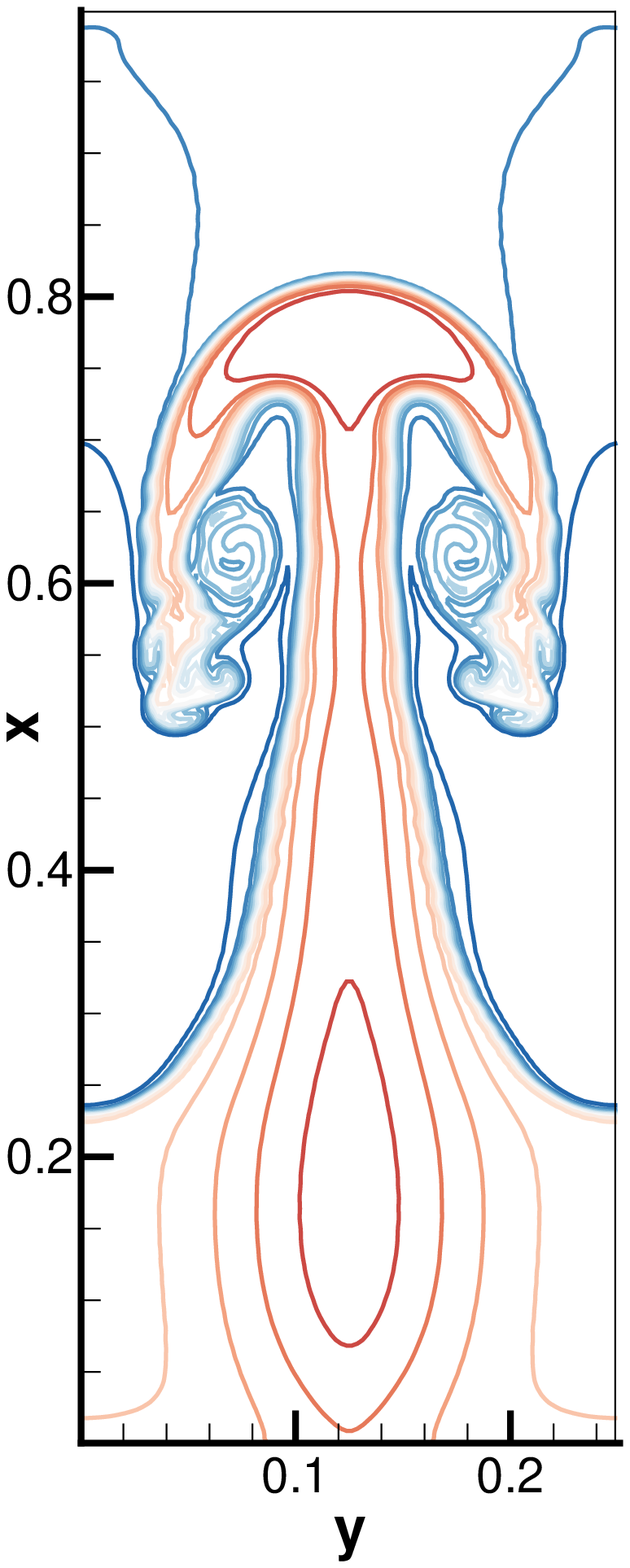}
\includegraphics[width=0.25\textwidth]{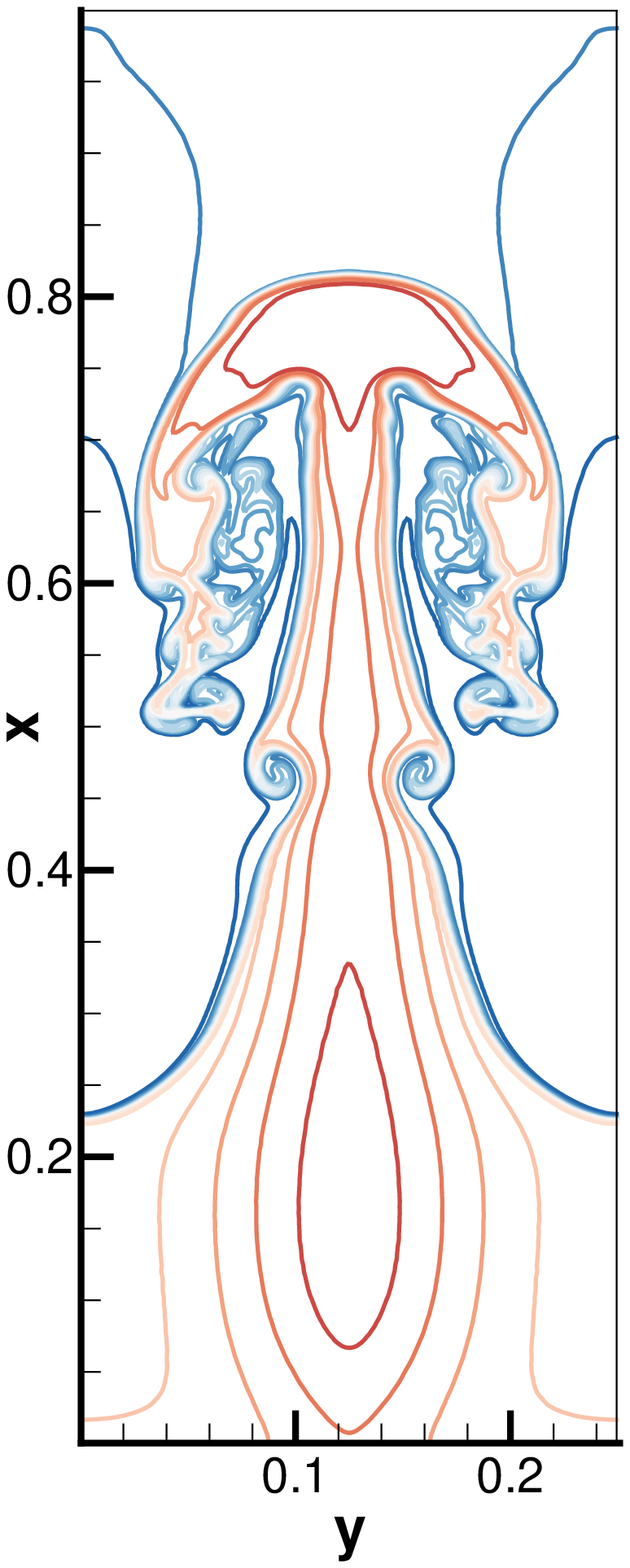}
\includegraphics[width=0.25\textwidth]{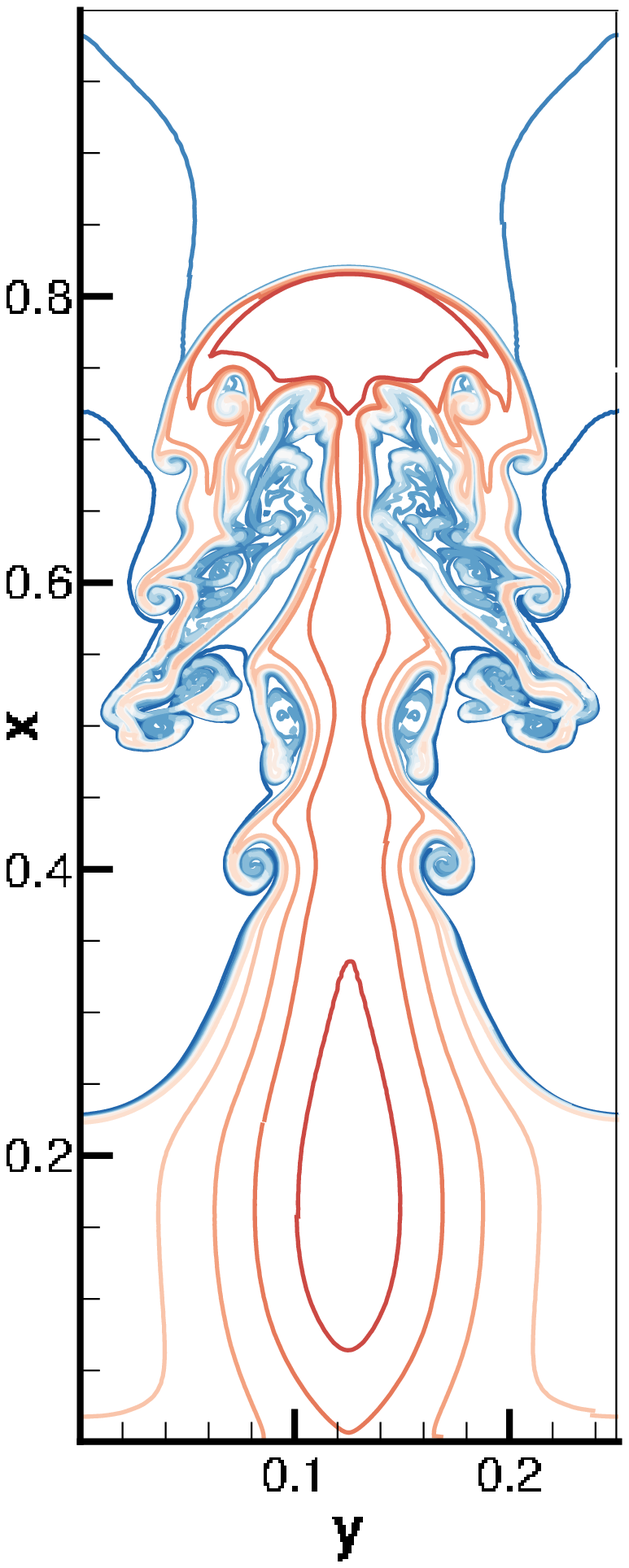}
\caption{\label{Rayleigh-Taylor} Rayleigh-Taylor instability: the density distribution with the mesh size $\Delta x=\Delta y=1/400, 1/800$ and $1/1600$ at $t=2.25$.}
\end{figure}

\subsection{Rayleigh-Taylor instability}
The last case is the Rayleigh-Taylor instability to test the performance of the scheme for the conservation laws with source terms. Rayleigh-Taylor instability happens on an interface between fluids with different densities when
an acceleration is directed from the heavy fluid to the light fluid. The instability has a fingering nature bubbles of light fluid rising into the ambient heavy fluid and spikes of heavy fluid falling into the light fluid. The initial condition of this problem \cite{Case-rt-instability} is given as follows
\begin{align*}
\begin{cases}
(\rho, U, V, p)=(2, 0, -0.025c\cos(8\pi x), 2y+1), x\leq 0.5,\\
(\rho, U, V, p)=(1, 0, -0.025c\cos(8\pi x), y+3/2), x> 0.5,
\end{cases}
\end{align*}
where $c=\sqrt{\displaystyle\frac{\gamma p}{\rho}}$ is the sound speed and $\gamma=5/3$. The computational domain is $[0,0.25]\times[0,1]$. The reflective boundary conditions are imposed for the left and right boundaries; at the top boundary, the flow values are set as $(\rho, U, V, p)=(1, 0, 0, 2.5)$, and at the bottom boundary, they are $(\rho, U, V, p)=(2, 0, 0, 1)$.  The source terms for the governing equations are $S(\textbf{u})=(0,0,\rho,\rho V)$.  The uniform meshes with $\Delta x=\Delta y=1/400, 1/800$ and $1/1600$ are used in the computation. The density distributions for the sixth-order scheme with WENO-JS reconstruction at $t=2.25$ are presented in Fig.\ref{Rayleigh-Taylor}. With the mesh refinement, the flow structures for the complicated flows are observed. It hints that current scheme may be suitable for the flow with interface instabilities as well.

\section{Conclusion}
In this paper, a new sixth-order WENO scheme was developed in finite volume framework for hyperbolic conservation laws. A new approach for selecting substencils and reconstruction procedure for the convex combination of candidate polynomials are introduced. Comparing with the classical WENO scheme, one order of improvement for the accuracy and better resolution in the region with discontinuity can be achieved with the same stencil by the current approach.
In this paper, only the sixth-order recursive WENO scheme is presented. Actually, as a WENO reconstruction approach of the complementary version of the classical WENO reconstruction, any even order accuracy schemes can be obtained based on the recursive WENO methodology. Numerical tests from the accuracy test to hypersonic flows are presented to validate accuracy and robustness of the recursive WENO scheme.

\section*{Acknowledgements}
The work of L. Pan is supported by China Postdoctoral Science Foundation (2016M600065).  The work of S. H. Wang is supported by NSAF (U1630247) and  NSFC (915303108).

\end{document}